\documentclass[11pt, a4paper]{article}
\usepackage{amsthm}
\usepackage{enumitem}
\usepackage{caption}
\usepackage{href-ul}
\usepackage{amsmath}
\usepackage{amssymb}
\usepackage{amscd}
\usepackage{graphicx}
\usepackage{geometry}
\usepackage{url,href-ul,cite}
\usepackage[noend]{algpseudocode}
\usepackage{algorithmicx,algorithm}
\usepackage{tikz-cd}
\usepackage{tikz}
\usepackage{authblk}
\newtheorem{thm}{Theorem}
\newtheorem{definition}{Definition}

\newtheorem{prop}{Proposition}

\newtheorem{example}{Example}

\title{IntComplex for high-order interactions} % The article title

\author[1]{Xiang Liu}
\author[1,2]{Ran Liu}
\author[1]{Jingyan Li}
\author[1]{Rongling Wu}
\author[1*]{Jie Wu}

\affil[1]{Beijing Institute of Mathematical Sciences and Applications}
\affil[2]{School of Mathematical Sciences, Beihang University}
\affil[*]{Corresponding author}

\date{} % Add a date here if you would like one to appear underneath the title block, use \today for the current date, leave empty for no date

\begin{document}
	
\maketitle % Print the title

\begin{abstract}
    Graphs serve as powerful tools for modeling pairwise interactions in diverse fields such as biology, material science, and social networks. However, they inherently overlook interactions involving more than two entities. Simplicial complexes and hypergraphs have emerged as prominent frameworks for modeling many-body interactions; nevertheless, they exhibit limitations in capturing specific high-order interactions, particularly those involving transitions from $n$-interactions to $m$-interactions. Addressing this gap, we propose IntComplex as an innovative framework to characterize such high-order interactions comprehensively. Our framework leverages homology theory to provide a quantitative representation of the topological structure inherent in such interactions. IntComplex is defined as a collection of interactions, each of which can be equivalently represented by a binary tree. Drawing inspiration from GLMY homology, we introduce homology for the detailed analysis of structural patterns formed by interactions across adjacent dimensions, $p$-layer homology to elucidate loop structures within $p$-interactions in specific dimensions, and multilayer homology to analyze loop structures of interactions across multiple dimensions. Furthermore, we introduce persistent homology through a filtration process and establish its stability to ensure robust quantitative analysis of these complex interactions. The proposed IntComplex framework establishes a foundational paradigm for the analysis of topological properties in high-order interactions, presenting significant potential to drive forward the advancements in the domain of complex network analysis.
\end{abstract}

%\tableofcontents
%\newpage

\section{Introduction}
A graph is a collection of vertices and the edges connecting them. In the context of modeling and analyzing systems, vertices correspond to the entities within the system, such as atoms in a protein or individuals in a social network. The edges denote the interactions between these entities, such as covalent bonds between atoms or social relationships between individuals. Graph theory has been successfully utilized in various fields, including biology, material sciences, social network analysis and computer sciences \cite{deo2017graph,gross2018graph,ortega2018graph,wang2019coevolution}. However, graphs implicitly ignore the high-order interactions since the edges only reflect the pairwise interactions.

Increasing evidence indicates that a single vertex can be influenced by multiple other vertices in a nonlinear manner, with such high-order interactions defying decomposition into simple pairwise interactions \cite{lambiotte2019networks,battiston2020networks,torres2021and,battiston2021physics} and the presence of high-order interactions may significantly impact the dynamics of the complex systems \cite{skardal2019abrupt,iacopini2019simplicial,neuhauser2020multibody,alvarez2021evolutionary}. Two prevalent frameworks for modeling these high-order interactions are simplicial complexes \cite{munkres2018elements} and hypergraphs \cite{berge1984hypergraphs}. A simplicial complex extends the concept of a graph by incorporating triangles, tetrahedrons, and their higher-dimensional counterparts, whereas a hypergraph generalizes a simplicial complex by allowing the presence of missing faces. Notably, by integrating homology theory from algebraic topology, these high-order frameworks based models have achieved significant success in various aspects of drug design \cite{ xia2014persistent,cang2018representability,liu2021hypergraph,liu2022dowker}, as well as in material design and discovery\cite{hiraoka2016hierarchical,lee2017quantifying, saadatfar2017pore,anand2022topological}.

\begin{figure}[h]
	\centering
	\includegraphics[width=0.95\textwidth]{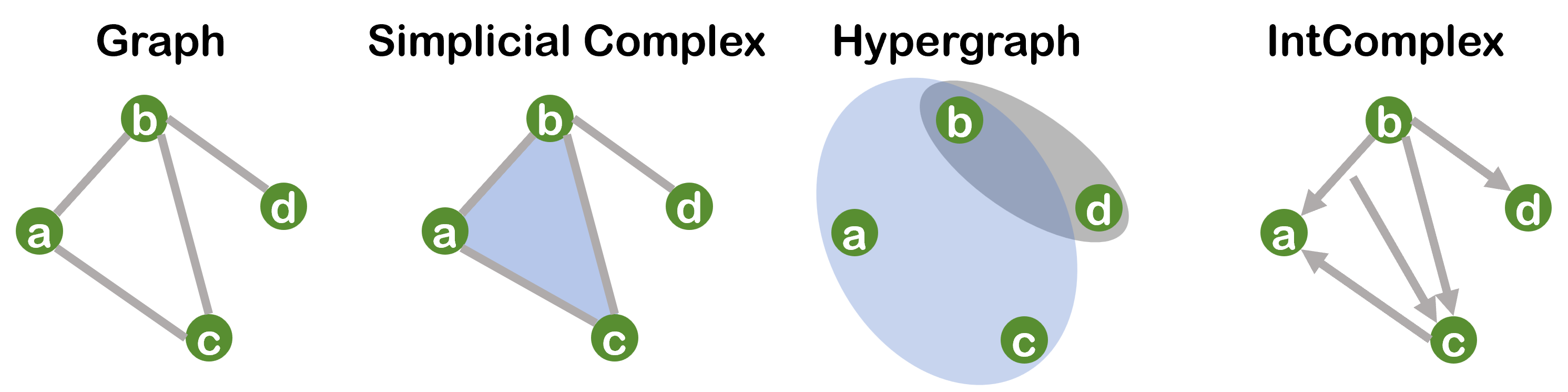}
	\caption{Comparison among a graph, simplicial complex, hypergraph and IntComplex. The simplicial complex and hypergraph can represent the high-order (3-order) interactions among {\bf a, b} and {\bf c} that cannot be captured in the graph representation. Specifically, the 3-order interaction among {\bf a, b} and {\bf c} are represented by a 2-simplex \{\textbf{a,b,c}\} and a 2-hyperedge \{\textbf{a,b,c}\} respectively. So the interaction from the interaction {\bf a} and {\bf b} to {\bf c} cannot be captured in the simplicial complex and hypergraph but can be represented as a 3-interaction (({\bf b}, {\bf a}), {\bf c}) in the IntComplex. }
	\label{fig:comparison}
\end{figure}

However, there are still high-order interactions that simplicial complexes and hypergraphs cannot model. In both simplicial complex and hypergraph representations, an $n$-body interaction is represented by a set of $n$ vertices, meaning there are no internal structures within these $n$-interactions. In practice, it is common for a pairwise interaction from A to B to interact with C, akin to how the relationship between a father and mother in a family naturally affects their son. Figure \ref{fig:comparison} shows the difference among graphs, simplicial complexes, hypergraphs and IntComplexes (We propose in this work).
As shown in Figure \ref{fig:comparison}, the 3-interaction from the 2-interaction of {\bf b} and {\bf a} to {\bf c} cannot be represented by the graph, simplicial complex and hypergraph but can be naturally modeled as a 3-interaction (({\bf b}, {\bf a}), {\bf c}) in the IntComplex. Consequently, simplicial complexes and hypergraphs are inadequate for modeling such high-order interactions. These high-order interactions are prevalent in various complex systems but have been explored in only a few studies\cite{wang2009genome,kenett2015partial,zhao2016part} and remains an important research area \cite{baptista2024mining}. And for the general high-order interaction from an $n$-interaction to an $m$-interaction, there is no mathematically well-defined modeling method as far as we know. 

Here we present IntComplex for the modeling of high-order interactions, an IntComplex is a collection of interactions and the interaction from an $n$-interactoin to an $m$-interaction is modeled as an $(m+n)$-interaction. Further, we give the homology theory for analyzing the topological structures within the IntComplex, including the standard homology for the structures between adjacent-order interactions, $p$-layer homology for the loop structures within $p$-interactions, and multilayer homology for the loop structures among interactions of several orders. By considering a filtration of the IntComplex, we propose the persistent homology of IntComplex to give a multiscale quantitative characterization of the topological structures within the IntComplex, and give the stability theorem of the persistent homology to ensure its robustness. The paper is organized as follows: in section \ref{section:intcomplex}, we give the construction of IntComplexes and two representations of an interaction and its face operation, in section \ref{section:homology}, we construct the homology theory of IntComplex and show that the homology is a simple homotopy invariant and a functor, in section \ref{section:persistent-homology}, we give the persistent homology of IntComplexes and its stability, in section \ref{section:experiment}
, we give some experimental tests, a conclusion is in section \ref{section:conclusion}.

\section{IntComplex}\label{section:intcomplex}
In this section, we present the definition of IntComplex along with two different representations of the interactions and the face operation of the interactions.
%\begin{definition}[Binary tree of a set]
%	Given a set $V$, a binary tree is called a p-binary-tree of V if the binary tree has exactly $p$ leaf nodes and all leaf nodes are from $V$, in particular, any element $v\in V$ is called a 1-binary tree of V.
%\end{definition}
\subsection{Interaction and IntComplex}
\begin{definition}[Interaction of a set]\label{defn:1}
	Given a set V,  $\forall v\in V$ is called an 1-interaction of V. Take any $s$-interaction $\sigma_s$ and $t$-interaction $\sigma_t$ of V, the ordered pair $(\sigma_s,\sigma_t)$ is called an $(t+s)$-interaction of V $(t,s>0)$. Denote the set of all $p$-interactions of V by $Int_p(V)(p>0)$, $Int(V)=\bigcup_pInt_p(V)$.
\end{definition}
From definition \ref{defn:1}, $\forall \sigma_p\in Int_p(V)(p>1)$, there exists unique $\sigma_i\in Int_i(V),\sigma_j\in Int_j(V)(i+j=p)$ such that $\sigma_p=(\sigma_i,\sigma_j).$ We call $\sigma_i$ and $\sigma_j$ the left daughter and right daughter of $\sigma_p$. We let the left and right daughters of any 1-interaction $\sigma_1$ be the empty set. Any interaction represents an interaction from its left daughter  to its right daughter. For example, $\sigma_p=(\sigma_i,\sigma_j)$ represents an interaction from $\sigma_i$ to $\sigma_j$.

\begin{definition}[IntComplex]
	An IntComplex $\mathcal{I}$ over a set $V$ is a nonempty subcollection of $Int(V)$. 
	%satisfying the property: $\forall \sigma\in \mathcal{I}$, the left daughter and right daughter of $\sigma$ are also in $\mathcal{I}$.
\end{definition}

Denote the set of $p$-interactions of $\mathcal{I}$ by $\mathcal{I}_p$. We also call $\mathcal{I}_p$ the $p$-th layer of $\mathcal{I}$. The elements of $\mathcal{I}_1$ are called vertices of $\mathcal{I}$. Clearly, $\mathcal{I}_1$ is a subset of $V$, we remove from $V$ all the non-vertices so that $V=\mathcal{I}_1$. The IntComplex $\mathcal{I}'$ is a subcomplex of the IntComplex $\mathcal{I}$ if $\mathcal{I}'_p\subset\mathcal{I}_p(p\geqslant1)$.
%When an IntComplex $\mathcal{I}$ is fixed, all the interactions of $\mathcal{I}$ are called allowed while all the interactions not in $\mathcal{I}$ are called nonallowed. The dimension of any $p$-interactions $\sigma_p$ is $p$, denoted by $dim(\sigma_p)=p$. 
%Note that $(a_i,a_j)=((a_i,a_j))=(((a_i,a_j)))=(...(a_i,a_j)...)$. We call the representation with least brackets the reduced one. For example, $(a_i,a_j)$ is reduced while $((a_i,a_j))$ is not reduced. Throughout this work, we always consider the reduced representation. If the representation is not reduced, we implicitlly change it to its reduced one. 
\begin{example}\rm
	Figure \ref{fig:hypernet} {\bf A} shows an IntComplex $\mathcal{I}$ over $V=\{a,b,c,d\}$. $\mathcal{I}_1=\{a,b,c,d\}$, $\mathcal{I}_2=\{(b,a),(d,d)\}$, $\mathcal{I}_3=\{((b,a),c),(c,(d,d))\}$, $\mathcal{I}_6=\{(((b,a),c),(c,(d,d)))\}$, $\mathcal{I}_p=\emptyset(p\ne 1,2,3,6)$
	
\begin{figure}[h]
	\centering
	\includegraphics[width=0.95\textwidth]{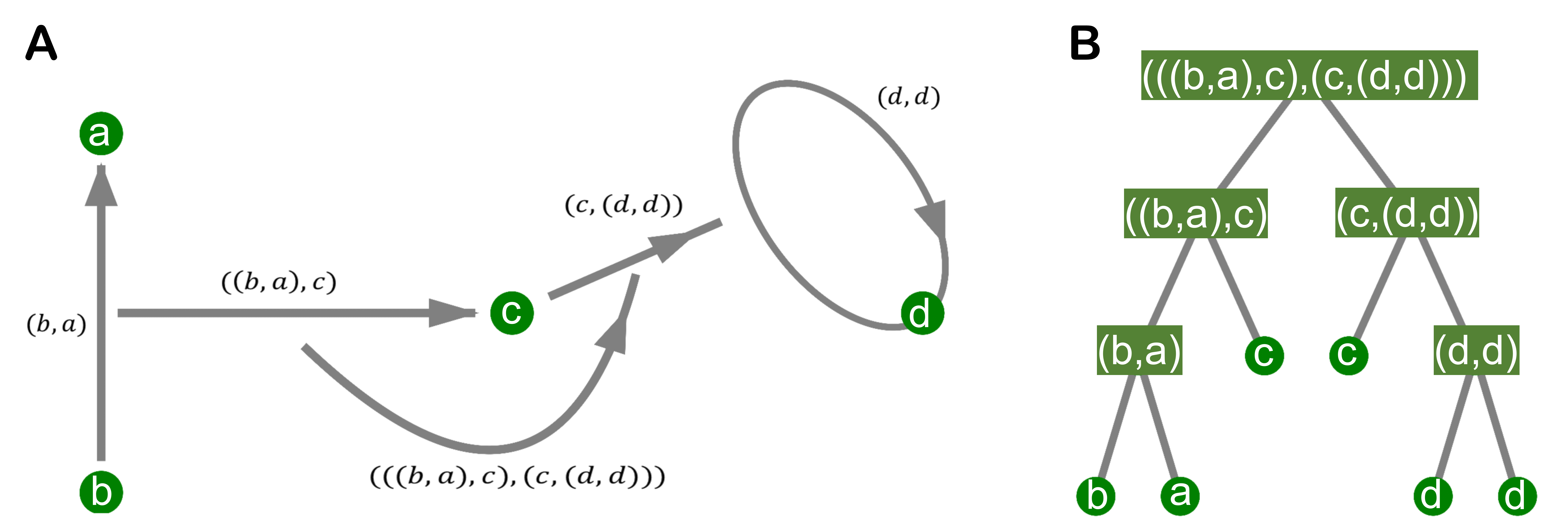}
	\caption{Illustration of the IntComplex and the binary tree representation of its interactions. {\bf A}: An IntComplex over $V=\{a,b,c,d\}$. {\bf B}: The binary tree representation of the interactions, every interaction $\sigma$ is represented by a binary tree whose root node is $\sigma$.}
	\label{fig:hypernet}
\end{figure}
\end{example}

\begin{example}\rm
	Given a directed graph $G$ with vertex set $V$ and edge set $E$. Let $\mathcal{I}_1=V,\mathcal{I}_2=E, \mathcal{I}_p=\emptyset(p>2)$, then $\mathcal{I}$ is an IntComplex. Hence any directed graph is a special IntComplex.
\end{example}
%\paragraph{Remark}
%\begin{enumerate}
		%\item Every element $a_n\in A_n$ can be represented by a sequence of $n$ elements of $A_1$ together with $"()"$ where  $"(a_i,a_j)"$ represents an interaction from $a_i$ to $a_j$.
		%\item Note that $((1,2),(3,4))$ is different from $(1,(2,(3,4)))$. $((1,2),(3,4))$ represents the interaction from 2-superarrow (1,2) to 2-superarrow (3,4). $(1,(2,(3,4)))$ represents the interaction from 1-superarrow 1 to 3-superarrow $(2,(3,4))$.
		%\item The element $a_n\in A_n$ also can be represented by a binary tree. For example\\
		%\includegraphics[width=0.8\textwidth]{ex3.png}
		%\item The tree representation and "()" representation of $e_n$ is uniquely determined each other. They are equivalent.
		%\item If $A_n=\emptyset~(n>1)$, G is a hypergraph.
		%\item If each $e_1\in A_1$ has only one verex and $A_n=\emptyset~(n>2)$, G is a digraph.
		%\item We can define homology for each $A_i(i\geqslant 1)$ separately. For example, $A_1$ is just ahypergraph, we can use the hypergraph homology. For $A_n$, let $A_n^i(i\geqslant n)$ represents the set of elements of $A_n$ consists of i vertices o f V. Then we have the graded set $A_n^n,A_n^{n+1},A_n^{n+2},...$, then the embedded homology can be defined. 
		%\item The homology for each $A_i$ can be seen as inner homology to characterize the inner structure within each $A_i$. The path homology of all $A_i$ can be seen as the global homology to characterize the global structure of the super-digraph.
%\end{enumerate}

\subsection{Interaction representation}
\subsubsection{Number Pair (NP) representation}
\begin{prop}\label{prop:1}
	Each n-interaction $\sigma_n\in Int_n(V)$ can be equivalently represented by a sequence of n ordered vertices $v_{1},v_{2},...,v_{n}$ and n-1 pairs of intergers $\{(l_i,r_i)|i=1,2,...,n-1\}$. And $\forall (l_{i_1},r_{i_1}),(l_{i_2},r_{i_2})$, they satisfy one of the followings
	\begin{enumerate}
		\item $[l_{i_1},r_{i_1}]\subsetneqq [l_{i_2},r_{i_2}]$ or $[l_{i_2},r_{i_2}]\subsetneqq [l_{i_1},r_{i_1}]$
		\item $[l_{i_1},r_{i_1}]\cap[l_{i_2},r_{i_2}]=r_{i_1}=l_{i_2}$ or  $[l_{i_1},r_{i_1}]\cap[l_{i_2},r_{i_2}]=l_{i_1}=r_{i_2}$ or $[l_{i_1},r_{i_1}]\cap[l_{i_2},r_{i_2}]=\emptyset$
	\end{enumerate}
	\begin{proof}
		By induction, we can prove that each $n$-interaction $\sigma_n$ can be uniquely represented by a sequence of $n$ ordered vertices $v_{1},v_{2},...,v_{n}$ and $n$-1 round brackets ``()". These $n$ ordered vertices provide $n$+1 positions for the brackets to be put in. We number these positions by 1, 2,..., $n$+1 from left to right. Then each bracket ``()'' uniquely corresponds to a pair of numbers $(i,j)$ that are the positions of the bracket.
		
		For the two properties of the number pairs. Each number pair corresponds to a bracket $``()"$, so the result follows by definition \ref{defn:1}. 
	\end{proof}
\end{prop}
We denote the $n$-interaction by $\sigma_n=[v_{1}...v_{n},\{(l_i,r_i)\}_{i=1}^{n-1}]$ and call it the NP representation of $\sigma_n$. The number pairs $\{(l_i,r_i)\}_{i=1}^{n-1}$ are positions of the $n-1$ brackets. For example, $[ba,\{(1,3)\}]$, $[dd,\{(1,3)\}]$, $[bac,\{(1,3),(1,4)\}]$, $[cdd,\{(2,4),(1,4)\}]$ and $[baccdd,\{(1,3),(1,4),(5,7),(4,7),(1,7)\}]$ are the NP representations of $(b,a)$, $(d,d)$, $((b,a),c)$, $(c,(d,d))$ and $(((b,a),c),(c,(d,d)))$ respectively.

\begin{algorithm}[h] 
	\caption{Int2Tree }
	\label{alg:bracket2tree}
	\hspace*{0.02in} {\bf Input:}
	an interaction $\sigma$ and the corresponding node $v$\\
	\hspace*{0.02in} {\bf Output:}
	the binary tree of $\sigma$
	\begin{algorithmic}[1]
		\If{$\sigma$ is an 1-interaction}
		\State return $v$
		\Else
		\State $\sigma=(\sigma_i,\sigma_j)$
		\State add two nodes $v_i,v_j$ as the left and right daughters of $v$
		\State return $T(v,{\rm Int2Tree}(\sigma_i,v_i),{\rm Int2Tree}(\sigma_j,v_j))$
  \EndIf
	\end{algorithmic}
\end{algorithm}

\begin{algorithm}[h] 
	\caption{Tree2Int }
	\label{alg:tree2bracket}
	\hspace*{0.02in} {\bf Input:}
	a binary tree $bt$\\
	\hspace*{0.02in} {\bf Output:}
	the interaction of $bt$
	\begin{algorithmic}[1]
		\If{$bt$ has only one node $v$ }
		%\State insert $v$ into the bracket to get ``$(v)$''
		\State return $v$
		\Else
		\State denote the left daughter and right daughter of $bt$ by $bt_1$ and $bt_2$ respectively.
		\State return ( Tree2Int($bt_1$), Tree2Int($bt_2$) )
		\EndIf
	\end{algorithmic}
\end{algorithm}
   
%\begin{definition}[Binary tree]
%	Given a set V, the binary tree over V is defined as:
%	\begin{itemize}
%		\item $\forall v\in V$, $v$ is a binary tree.
%		\item if $v_1,v_2$ are binary trees, then the tree by adding a root $r$ connected to the left to $v_1$ and to the right to $v_2$ by edges is a binary tree, denoted by $(v_1,r,v_2)$.
%	\end{itemize}
%	$v_1,v_2$ are called the left daughter and right daughter of $(v_1,r,v_2)$ respectively.
%\end{definition}
\subsubsection{Binary Tree (BT) representation}
A binary tree is a tree in which each node has at most two daughter nodes. We denote a binary tree with $v$ as root node and $v_i,v_j$ as the left and right daughters by $T(v,v_i,v_j)$.
Each $n$-interaction $\sigma_n\in Int_n(V)$ can be represented as a binary tree by Algorithm \ref{alg:bracket2tree}. %From Proposition \ref{prop:1} we know, $\sigma_n$ can be uniquely represented as n ordered vertices $v_{1},v_{2},...,v_{n}$. These vertices are one-to-one correspond to the leavies of the binary tree. 
On the contrary, every binary tree $bt$ can be represented by an interaction by Algorithm \ref{alg:tree2bracket}.	 
Let $T(\sigma_n)=\rm Int2Tree(\sigma_n)$, $S(bt)=\rm Tree2Int(bt)$. 

Figure \ref{fig:hypernet} {\bf B} gives an example of the binary tree representation of the interactions in {\bf A}. Each interaction $\sigma$ is represented by a binary tree whose root node is $\sigma$. For example, the 2-interaction $(b,a)$ is represented by the binary tree with two leaf nodes $b$ and $a$. The 3-interaction $(c,(d,d))$ is represented by the binary tree with three leaf nodes $c$, $d$ and $d$. 

\begin{prop}\label{prop:9}
	For any $n$-interaction $\sigma_n$, $\sigma_n=S( T(\sigma_n) )$, for any binary tree $bt$, $bt=T( S(bt) )$.
	\begin{proof}
		Algorithm \ref{alg:bracket2tree} and Algorithm \ref{alg:tree2bracket} are inverse to each other.
	\end{proof}
\end{prop}
From proposition \ref{prop:9}, we know each $n$-interaction has an equivalent NP representation containing $n$ vertices and $n-1$ pairs of position numbers,  and an equivalent BT representation with $n$ leaf nodes. The binary tree representaiton can be derived by Algorithm \ref{alg:bracket2tree}.  We call $T(\sigma_n)$ the binary tree of $\sigma_n$.

\subsection{Face}
Given an interaction $\sigma$, consider its binary tree $T(\sigma)$, the $i$-th (from left to right) leaf node $v_i$ of $T(\sigma)$ has a minimal binary tree that contains $v_i$ as a daughter, denoted by $T_i$. Removing the $i$-th leaf node of $T(\sigma)$ is equivalent to repalce this subtree $T_i$ by the other daughter of $T_i$. In the NP representation, it is to remove the vertex $v_i$ and the number pair corrsponding to $T_i$. 

\begin{prop}\label{prop:3}
	Consider an $n$-interaction $\sigma_n=[v_{1}...v_{n},\{(l_i,r_i)\}_{i=1}^{n-1}]$. $\forall v_{p}(1\leqslant p\leqslant n)$,  $\exists q$ such that $l_q\leqslant p< r_q$. Among such numbers q, there is a unique $q'$ such that $r_{q'}-l_{q'}$ is the smallest one.
	\begin{proof}
		The existance is obvious. For the uniqueness, suppose there are two number pairs $(l_{q_1},r_{q_1})$, $(l_{q_2},r_{q_2})$ such that $l_{q_1}\leqslant p< r_{q_1}$, $l_{q_2}\leqslant p< r_{q_2}$, then the interaction of these two pairs are not just one number, which contradicts with Proposition \ref{prop:1}.
	\end{proof} 
\end{prop}
We denote the number $q'$ by $T(v_{i},\sigma_n)=q'$. The unique number pair of $v_i$ just corresponds to the minimal binary tree $T_i$ of $v_i$. Now we can give the face operation.

\begin{definition}[Face]
	Given an n-interaction $\sigma_n=[v_1...v_n,\{(l_i,r_i)\}_{i=1}^{n-1}]$ over V, let $M_{\sigma_n}^j=\{(l_i,r_i)|1\leqslant i\leqslant n-1, i\ne T(v_j,\sigma_n)\}$. For each $(l_t,r_t)\in M^j_{\sigma_n}$, if $l_t>j$, replace $l_t$ by $l_t-1$, if $r_t>j$, replace $r_t$ by $r_t-1$. Perform this replacement for all elements of $M^j_{\sigma_n}$, denote the resulting set by $N^j_{\sigma_n}$. Then $[v_1...v_{j-1}v_{j+1}...v_n,N_{\sigma_n}^j]$ is an $(n-1)$-interaction, define it as the $j$-th face of $\sigma_n$, denoted by $F_j(\sigma_n)$. Given a binary tree T, define the $j$-th face of T to be the binary tree derived from T by removing the $j$-th leaf node, denoted by $F_j(T)$.
\end{definition}
The $n$-interaction $\sigma_n$ has $n$ faces $\{ F_1(\sigma_n),F_2(\sigma_n),...,F_n(\sigma_n)\}$. The binary tree of $j$-th face of $\sigma_n$ is exactly the binary tree derived from $T(\sigma_n)$ by removing the $j$-th leaf node. That is, $F_j(T(\sigma_n))=T(F_j(\sigma_n)).$

\begin{prop}\label{prop:11}
	$\forall i<j,F_iF_j=F_{j-1}F_i$.
	\begin{proof}
		For any binary tree $bt$, $F_iF_j(bt)$ is the binary tree from $bt$ by first removing the $j$-th and then the $i$-th leaf nodes. $F_{j-1}F_i(bt)$ is the binary tree from $bt$ by first removing the $i$-th and then the $(j-1)$-th leaf nodes. Hence $F_iF_j(bt)=F_{j-1}F_i(bt)$.
		
		$\forall \sigma_n\in Int_n(V)$, We have $F_iF_j(T(\sigma_n))=F_{j-1}F_iT(\sigma_n)$, so $S(F_iF_jT(\sigma_n))=SF_{j-1}F_iT(\sigma_n)$, $STF_iF_j(\sigma_n)=STF_{j-1}F_i(\sigma_n)$, that is, $F_iF_j=F_{j-1}F_i$.
	%	$\forall \sigma_n=[v_1...v_n,\{l_t,r_t\}_{t=1}^{n-1}]\in Int_n(V)$, consider all cases of the positions of $v_i,v_j$ in the bracket representation.
	%	\begin{enumerate}
	%		\item $(L~((v_i,v_j),\square)~R)$, 
	%		\item $(L~(\square,(v_i,v_j))~R)$. 
	%		\item $(L~(v_i,(\square,v_j))~R)$, 
	%		\item $(L~(v_i,(v_j,\square))~R)$
	%		\item $(L~((v_i,\square),v_j)~R)$
	%		\item $(L~(v_i,\square)~M~(v_j,\square)~R)$,  
	%		\item $(L~(v_i,\square)~M~(\square,v_j)~R)$, 
	%		\item $(L~(\square,v_i)~M~(v_j,\square)~R)$,
	%		\item $(L~(\square,v_i)~M~(\square,v_j)~R)$,
	%	\end{enumerate}
	\end{proof}
\end{prop}

\section{Homology of IntComplex}\label{section:homology}
In this section, we introduce three types of homology groups for the IntComplex. First, we present homology for the structure composed of adjacent-layer interactions. Next, we describe layer-homology for structures formed by interactions within specific layers. Lastly, we discuss multilayer-homology for structures arising from interactions across multiple layers. Specifically, our focus lies on analyzing loop structures for layer-homology and multilayer-homology.
\subsection{Boundary operator and Homology}
Fix a field $\mathbb{F}$ as the coefficient. Given an IntComplex $\mathcal{I}$ over $V$, %$A_1=\{a_1^1,a_1^2,...,a_1^n\}$, Let $\Lambda_p(G)(1\leqslant p\leqslant n)$ represent the set of all binary trees that have $p$ leaf nodes and all $p$ leaf nodes are from $A_1$, any element $a_p\in \Lambda_p(G)$ also can be represented as $a_p=[a_1^1...a_1^p,\{(l_i^{a_p},r_i^{a_p})\}_{i=1}^{p-1}]$, we also call elements of $\Lambda_p(G)$ $p$-superarrows. We have $A_p\subset \Lambda_p(G)$. 
Let $\Lambda_p(\mathcal{I})$ be the $\mathbb{F}$-vector space spanned by $Int_{p}(V)(p>0)$.

\begin{figure}[h]
	\centering
	\includegraphics[width=0.8\textwidth]{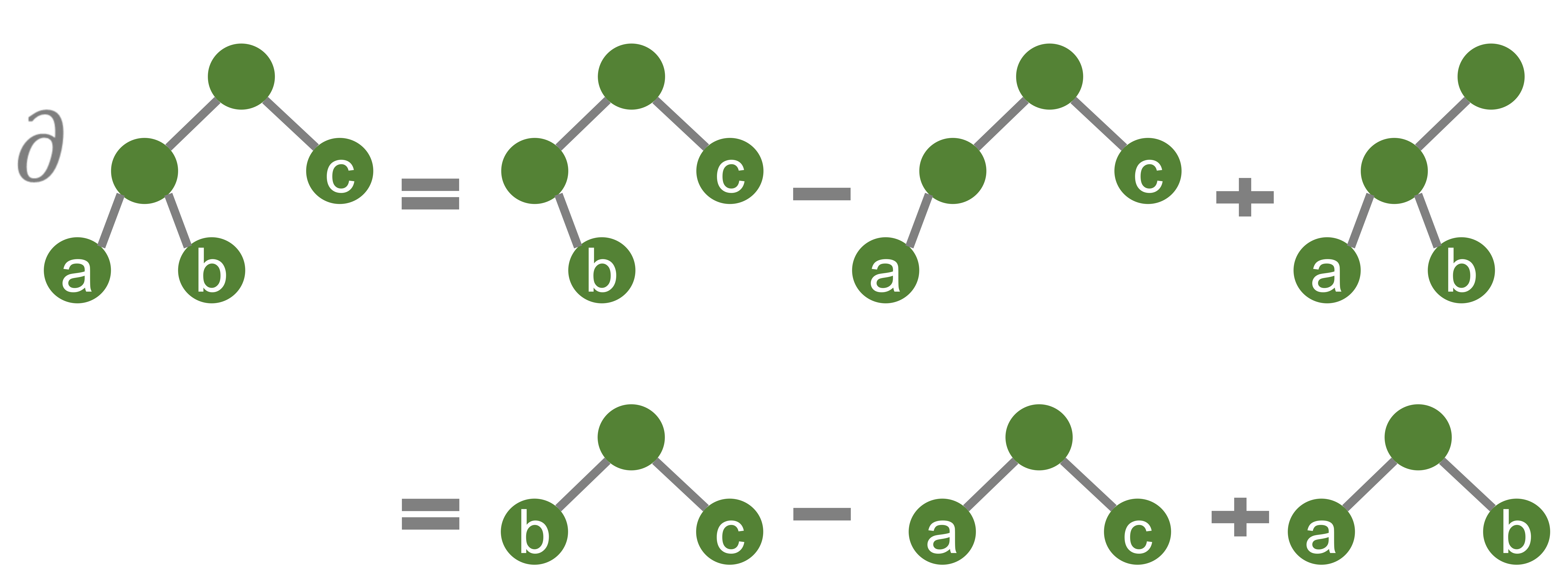}
	\caption{Illustration of the boundary operation bybinary tree.}
	\label{fig:boundary}
\end{figure}

%$\forall \sigma_k=[v_1v_2...v_{k},\{(l_i^{\sigma_k},r_i^{\sigma_k})\}_{i=1}^{k-1}]$, let $M_{\sigma_k}^i=\{(l_i^{\sigma_k},r_i^{\sigma_k})|1\leqslant i\leqslant k-1,i\ne T(v_i,\sigma_k)\}$. For each $(l_t^{\sigma_k},r_t^{\sigma_k})\in M^i_{a_k}$, if $l_t^{\sigma_k}>i$, replace $l_t^{\sigma_k}$ by $l_t^{\sigma_k}-1$, if $r_t^{\sigma_k}>i$, replace $r_t^{\sigma_k}$ by $r_t^{\sigma_k}-1$. Perform this replacement for all elements of $M^i_{\sigma_k}$, denote the resulting set by $N^i_{\sigma_k}$.
\begin{definition}[Boundary operator]
	For the $n$-interaction $\sigma_n=[v_1...v_n,\{(l_i,r_i)\}_{i=1}^{n-1}]\in Int_n(V)$, We define the boundary operator $\partial_n:\Lambda_n(\mathcal{I})\rightarrow \Lambda_{n-1}(\mathcal{I})$ as follows:
	$$\partial_n(\sigma_n)=\sum_{j=1}^{k}(-1)^{j+1}F_j(\sigma_n)$$
 That is,
	$$\partial_n([v_1...v_n,\{(l_i,r_i)\}_{i=1}^{n-1}])=\sum_{j=1}^k(-1)^{j+1}[v_1...v_{i-1}v_{i+1}...v_n,N_{\sigma_n}^j]$$
\end{definition}
Further, we define $\Lambda_{0}(\mathcal{I})=0$ and $\partial:\Lambda_1(\mathcal{I})\rightarrow\Lambda_{0}(\mathcal{I})$ to be zero.
In the binary tree representation, we have $$\partial_n(T(\sigma_n))=\sum_{j=1}^{n}(-1)^{j+1}T(F_j(\sigma_n))$$
where $T(F_j(\sigma_n))$ is derived by removing the $j$-th leaf node of $T(\sigma_n)$. 
Figure \ref{fig:boundary} shows the binary tree representation of the boundary operator for $\partial((a,b),c)=((b),c)-((a),c)+(a,b)=(b,c)-(a,c)+(a,b)$\\

\begin{prop}
	$\partial\partial=0$
	\begin{proof}
		$\forall \sigma_n\in Int_n(V)$,
		\begin{equation}\nonumber
			\begin{split}
				\partial(\partial(\sigma_n))&=\partial(\sum_{i=1}^n(-1)^iF_i(\sigma_n))\\
				&=\sum_{i=1}^n(-1)^i\partial(F_i(\sigma_n))\\
				&=\sum_{i=1}^n(-1)^i\sum_{j=1}^{i-1}(-1)^{j+1}F_j(F_i(\sigma_n))+\sum_{i=1}^n(-1)^i\sum_{j=i+1}^n(-1)^jF_{j-1}(F_i(\sigma_n))\\
				&=\sum_{1\leqslant j<i\leqslant n}(-1)^{i+j+1}F_j(F_i(\sigma_n))+\sum_{1\leqslant i<j\leqslant n}(-1)^{i+j}F_{j-1}(F_i(\sigma_n))\\
				&=\sum_{1\leqslant i<j\leqslant n}(-1)^{i+j+1}F_i(F_j(\sigma_n))+\sum_{1\leqslant i<j\leqslant n}(-1)^{i+j}F_{j-1}(F_i(\sigma_n))~(shift~i~and~j~of~the~first~part)\\
				&=\sum_{1\leqslant i<j\leqslant n}(-1)^{i+j+1}F_i(F_j(\sigma_n))-\sum_{1\leqslant i<j\leqslant n}(-1)^{i+j+1}F_{i}(F_j(\sigma_n))~(by ~proposition~\ref{prop:11})\\
				&=0
			\end{split}
		\end{equation}

	%$N_{a_k}^{i,j}$ represents the set of number pairs by firstly removing the smallest pair that contains $j$ from $N_{a_k}^i$ then replacing numbers $k$ larger than $j$ by $k-1$. Hence, it suffices to prove that for any $i\ne j$, $N_{a_k}^{i,j}=N_{a_k}^{j,i}$. In the binary tree representation, this refers to firstly removing $i$-th leaf node then removing $j$-th leaf node and firstly removing $j$-th leaf node then removing $i$-th leaf node will give the same binary tree, which is obvious.  		
	\end{proof}
\end{prop}
Consequently, we get a chain complex $$...\rightarrow \Lambda_{p+1}(\mathcal{I})\xrightarrow{\partial}\Lambda_p(\mathcal{I})\xrightarrow{\partial}\Lambda_{p-1}(\mathcal{I})\xrightarrow...\rightarrow\Lambda_1(\mathcal{I})\xrightarrow{0}\Lambda_0(\mathcal{I})=0$$
Denoted by $(\Lambda_*(\mathcal{I}),\partial_*)$. Let $\mathcal{A}_p(\mathcal{I})$ be the $\mathbb{F}$-vector space spanned by $\mathcal{I}_p$, $\mathcal{A}_p(\mathcal{I})\subset \Lambda_p(\mathcal{I})$, but generally, $\partial(\mathcal{A}_p(\mathcal{I}))\not\subset \mathcal{A}_{p-1}(\mathcal{I})$. Considering the following subspace of $\mathcal{A}_p(\mathcal{I})$
$$\Omega_p(\mathcal{I})=\{\alpha_p\in \mathcal{A}_p(\mathcal{I})|\partial(\alpha_p)\in \mathcal{A}_{p-1}(\mathcal{I})\}$$
We have $\partial(\Omega_p(\mathcal{I}))\subset \Omega_{p-1}(\mathcal{I})$, so we get a chain complex
$$...\rightarrow \Omega_{p+1}(\mathcal{I})\xrightarrow{\partial}\Omega_p(\mathcal{I})\xrightarrow{\partial}\Omega_{p-1}(\mathcal{I})\rightarrow...\rightarrow \Omega_1(\mathcal{I})\xrightarrow{0}0$$
Denoted by $(\Omega_*(\mathcal{I}),\partial_*(\mathcal{I}))$.
\begin{definition}[Homology of IntComplex]
	Define the p-th homology of the IntComplex $\mathcal{I}$ as the p-th homology of $(\Omega_*(\mathcal{I}),\partial_*(\mathcal{I}))$
	$$H_p(\mathcal{I})=H_p(\Omega_*(\mathcal{I}),\partial_*(\mathcal{I}))=Ker\partial_{p}|_{\Omega_*(\mathcal{I})}/Im\partial_{p+1}|_{\Omega_*(\mathcal{I})}$$
	Define the $p$-th betti number of $\mathcal{I}$ as the dimension of $H_p(\mathcal{I})$, denoted by $\beta_p=dim(H_p(\mathcal{I}))$.
\end{definition}

\subsection{Simple Homotopy type}
\begin{definition}
    Given an IntComplex $\mathcal{I}$, a pair of interactions $(\sigma_p,\tau_{p+1})$ is called free if $\sigma_p$ is a face of $\tau_{p+1}$ and not face of any other $(p+1)$-interactions. 
\end{definition}

\begin{definition}
    Given an IntComplex $\mathcal{I}$, assume $(\sigma,\tau)$ is a free pair of $\mathcal{I}$, $\mathcal{I}'=\mathcal{I}\setminus\{\sigma,\tau\}$ is also an IntComplex. We call $\mathcal{I}'$ is derived by an collapse from $\mathcal{I}$ and $\mathcal{I}$ is derived from $\mathcal{I}'$ by an expansion.
\end{definition}

\begin{definition}
    Given two IntComplexes $\mathcal{I}$ and $\mathcal{I}'$, we say they have the same simple homotopy type if $\mathcal{I}$ can be derived from $\mathcal{I}'$ by finite collapses and expansions.
\end{definition}

\begin{thm}
    If two IntComplexes $\mathcal{I}$ and $\mathcal{I}'$ have the same simple homotopy type, then for any $ p\geqslant0$
    $$H_p(\mathcal{I})\cong H_p(\mathcal{I}')$$
    \begin{proof}
        Assume $\mathcal{I}'$ is derived from $\mathcal{I}$ by a collapse of the free pair $(\sigma_n,\tau_{n+1})$ where $\sigma_n$ is an $n$-interaction and $\tau$ is an $(n+1)$-interaction. It suffices to prove that $H_p(\mathcal{I})\cong H_p(\mathcal{I}')$.
        We have the chain complexes $(\Omega_*(\mathcal{I},\partial_*),(\Omega_*(\mathcal{I}',\partial_*)$ for $\mathcal{I}$ and $\mathcal{I}'$ respectively. Note that $(\Omega_*(\mathcal{I}'),\partial_*)$ is a subchain complex of $(\Omega_*(\mathcal{I}),\partial_*)$, so we have the quotient complex $(\Omega_*(\mathcal{I}/\mathcal{I}'),\partial_*)$ where
        \begin{equation}\nonumber
            \Omega_p(\mathcal{I}/\mathcal{I}')=\Omega_p(\mathcal{I})/\Omega_p(\mathcal{I}')
        \end{equation}
        and the exact sequence of chain complex
        \begin{equation}\nonumber
            0\rightarrow \Omega_*(\mathcal{I}')\xrightarrow{i}\Omega_*(\mathcal{I})\xrightarrow{j}\Omega_*(\mathcal{I}/\mathcal{I}')\rightarrow 0
        \end{equation}
        where $i$ is the subcomplex injection and $j$ is the quotient projection. Consequently, it suffices to prove the quotient complex $(\Omega_*(\mathcal{I}/\mathcal{I}'),\partial_*)$ is acyclic, that is, $H_p(\Omega_*(\mathcal{I}/\mathcal{I}'))=0(p\geqslant 0)$. 
        \begin{enumerate}
            \item if there not exists $\alpha\in \Omega_n(\mathcal{I})$ such that $\sigma_n\in \alpha$, $\tau_{n+1}\notin \Omega_{n+1}(\mathcal{I})$. In this case, $\Omega_p(\mathcal{I})=\Omega_p(\mathcal{I}')(\forall p\geqslant 0)$, so $\Omega_*(\mathcal{I}/\mathcal{I}')=0(\forall p\geqslant 0)$, it follows that the quotient complex is acyclic.
            \item if there exists $\alpha\in \Omega_n(\mathcal{I})$ such that $\sigma_n\in\alpha$, there is $\gamma\in \Omega_{n+1}(\mathcal{I})$ such that $\tau_{n+1}\in\gamma$ and $\partial(\gamma)=\alpha$. In this case, 
            \begin{equation}\nonumber
                \begin{cases}
                \Omega_p(\mathcal{I})=\Omega_p(\mathcal{I}')& p<n\\
                \Omega_p(\mathcal{I})=\Omega_p(\mathcal{I}')\oplus<\alpha>& p=n\\
                \Omega_p(\mathcal{I})=\Omega_p(\mathcal{I}')\oplus<\gamma> &p=n+1\\
                \Omega_p(\mathcal{I})=\Omega_p(\mathcal{I}')& p>n+1\\
                
                \end{cases}
            \end{equation}
            So the quotient complex has nonzero chain group only in dimension $n$ and $n+1$. We have the following quotient complex
            \begin{equation}\nonumber
                0\rightarrow<\gamma>\xrightarrow{\partial}<\alpha>\rightarrow0
            \end{equation}
            $Ker\partial=0$, so $H_{n+1}(\Omega_*(\mathcal{I}/\mathcal{I}'))=0$. $\partial(\gamma)=\alpha$, so $H_{n}(\mathcal{I}/\mathcal{I}')=0$. It follows that he quotient complex is acyclic.
        \end{enumerate}
    \end{proof}
\end{thm}

\subsection{Connected components and $H_1$}
Given an IntComplex $\mathcal{I}$, a path of $\mathcal{I}$ is a sequence of vertices $v_0,v_1,...,v_n$ such that $(v_i,v_{i+1})\in \mathcal{I}_2 $ or $(v_{i+1},v_i)\in \mathcal{I}_2$ for all $0\leqslant i<n$.  Two vertices $u,v$ are called in the same connected component if there is a path connecting $u,v$. $\mathcal{I}$ is called connected if it only has one connected component. 
\begin{definition}[Disjoint union]
	Given two IntComplexs $\mathcal{I}'$ and $\mathcal{I}''$ over $V'$ and $V''$. The disjoint union of $\mathcal{I}$ and $\mathcal{I}'$ is the IntComplex $\mathcal{I}=\mathcal{I}'\coprod \mathcal{I}''$ with vertex set $V\coprod V'$, and $\mathcal{I}_n=\mathcal{I}_n'\coprod \mathcal{I}_n''$. Here $\coprod$ means disjoint union. 
\end{definition}
\begin{prop}\label{prop:10}
	Given two IntComplex $\mathcal{I}$ and $\mathcal{I}'$ over $V$ and $V'$, we have $$H_p(\mathcal{I}\coprod\mathcal{I}')=H_p(\mathcal{I})\oplus H_p(\mathcal{I}')$$
	\begin{proof}
		This is from the obvious identity: $\Omega_p(\mathcal{I}\coprod \mathcal{I}')=\Omega_p(\mathcal{I})\oplus\Omega_p(\mathcal{I}')$, $\partial(\Omega_p(\mathcal{I}))\subset \Omega_{p-1}(\mathcal{I})$, $\partial(\Omega_p(\mathcal{I}'))\subset \Omega_{p-1}(\mathcal{I}')$. 
	\end{proof}
\end{prop}
\begin{prop}\label{prop:12}
	For a connected IntComplex $\mathcal{I}$ over $V$, $dim(H_1(\mathcal{I}))=1$.
	\begin{proof}
		Assume $V=\{v_1,...,v_n\}$. $\forall v_i,v_j(i\ne j)$, there is a path connecting $v_i$ and $v_j$, which means $\exists \sigma_2\in \mathcal{I}_2$ such that $\partial(\sigma_2)=v_j-v_i$. So $v_j-v_i\in Im(\partial_2)(\forall i\ne j)$. So $dim(Ker(\partial_1))-dim(Im(\partial_2))=dim(\mathcal{I}_1)-dim(Im(\partial_2))\leqslant 1$. Note that $v_1\in \mathcal{I}_1,v_1\notin Im(\partial_2)$. So $dim(\mathcal{I}_1)-dim(Im(\partial_2))=1$. That is $dim(H_1(\mathcal{I}))=1$.
	\end{proof}
\end{prop}
\begin{prop}
	For any IntComplex $\mathcal{I}$, we have 
	$$dim(H_1(\mathcal{I}))=C$$
	where $C$ is the number of connected components of $\mathcal{I}$.
	\begin{proof}
		This is obvious from Proposition \ref{prop:10} and Proposotion \ref{prop:12}.
	\end{proof}
\end{prop}
\subsection{$p$-Layer-Homology}
\begin{definition}[Layer-homology]
	Given an IntComplex $\mathcal{I}$, consider $\mathcal{I}_p(p>1)$, let $V^\mathcal{I}_p=\{\tau|\tau$ is the daughter of $\sigma,\sigma\in \mathcal{I}_p\}$, $E^\mathcal{I}_p=\{\sigma=(\sigma_1,\sigma_2)|\sigma\in \mathcal{I}_p\}$, then $G^\mathcal{I}_p=(V^\mathcal{I}_p,E^\mathcal{I}_p)$ forms an IntComplex. We define the 1-th and 2-th homology of the IntComplex $G^\mathcal{I}_p$ as the  $p$-layer-homology of $\mathcal{I}$.
	$$H_p^L(\mathcal{I})=(H_1(G^\mathcal{I}_p),H_2(G^\mathcal{I}_p))$$
	denoted by $H_p^L(\mathcal{I})=(H_p^{L_1}(\mathcal{I}),H_p^{L_2}(\mathcal{I}))$. Define the $p$-layer betti number as the dimension of $p$-layer-homology, denoted by $\beta_p^L=(\beta_p^{L_1},\beta_p^{L_2})=(dim(H_p^{L_1}(\mathcal{I})),dim(H_p^{L_2}(\mathcal{I})))$
\end{definition}
%\begin{prop}
%	Given an IntComplex $\mathcal{I}$, we have,
%	$$dim(H_p^{L_1}(\mathcal{I}))-dim(H_p^{L_2}(\mathcal{I}))=|V_p^{\mathcal{I}}|-|E_p^\mathcal{I}|$$
%	where $|-|$ is the cardinal.  
%	\begin{proof}
%		This is a obvious fact from graph homology.
%	\end{proof}
%\end{prop}
The $p$-layer-homology of $\mathcal{I}$ reflects the structure of $p$-interactions of $\mathcal{I}$. $H_p^{L_1}(\mathcal{I})$ corresponds to the connected components clustered by the daughters of $p$-interactions, $H_p^{L_2}(\mathcal{I})$ correspond to the loop structures formed by $p$-interactions.

\begin{example}\label{exm:3}
	Consider the IntComplex in Figure \ref{fig:3}, $\mathcal{I}_1=V=\{1,2,3,4\}$,  $\mathcal{I}_2=\{(1,2),(2,3),(3,4),(1,4)\}$, $\mathcal{I}_p=\emptyset(p>2)$. The homology are shown in the right table.
	$\beta_1=1$ corresponds to the connected component $\{1,2,3,4\}$. $\beta_2=1$ corresponds to the loop formed by $\{(1,2),(2,3),(3,4),(1,4)\}$. For the layer-homology, $\beta_2^{L_1}=1$ corresponds to the connected component $\{1,2,3,4\}$. $\beta_2^{L_2}=1$ corresponds to the loop formed by the 2-interactions $\{(1,2),(2,3),(3,4),(1,4)\}$.

	Generally, for the IntComplex $\mathcal{I}$ over $V=\{1,2,...,n\}$ with $\mathcal{I}_2=\{[1,2],[2,3],...,[n-1,n],[1,n]\}$ where $[i,j]$ represents the 2-interaction $(i,j)$ or $(j,i)$. The homology are exactly same with the above example.
	\begin{figure}[ht]
		\centering
		\begin{minipage}{0.3\linewidth}
			\flushright
			\includegraphics[width=0.9\textwidth]{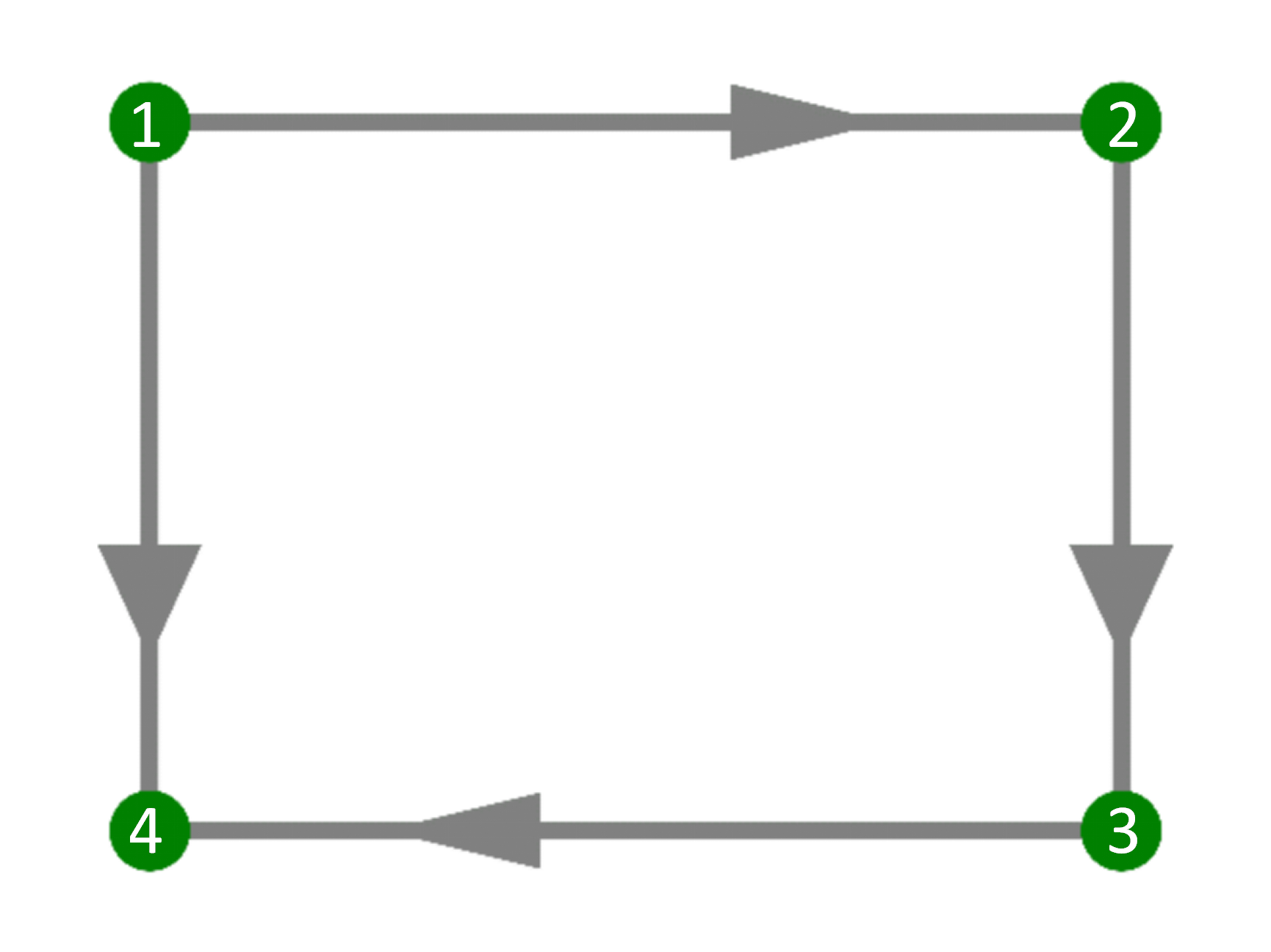}
			%\caption{An IntComplex.}
			%\label{fig:example1}
		\end{minipage}
		\hfill
		\begin{minipage}{0.6\linewidth}
			\centering
			\begin{tabular}{|c|c|c|c|}   
				\hline  &$\beta_1$ & $\beta_2$ & $\beta_p(p>2)$ \\   
				\hline
				homology   &  1 & 1 & 0 \\
				\hline
				Layer-Homology  &-& (1,1) & (0,0)\\
				\hline  	
			\end{tabular}   
			%\captionof{table}{Homology of the IntComplex in Figure \ref{fig:example1}. }
			%\label{table:example1}
		\end{minipage}
		\caption{An IntComplex and it homology.}
		\label{fig:3}
	\end{figure}
\end{example}

\begin{example}\label{exm:4}
	Consider the IntComplex in Figure \ref{fig:4}, $\mathcal{I}_1=V=\{1,2,3,4\}$, $\mathcal{I}_2=\{(1,2),(2,1)\}$, $\mathcal{I}_3=\{((1,2),3),((2,1),3),((1,2),4),((2,1),4)\}$, $\mathcal{I}_p=\emptyset(p>3)$. The homology are shown in the right table. $\beta_1=3$ corresponds to the three connected components $\{1,2\}, \{3\}$ and $\{4\}$. $\beta_3=1$ corresponds to the two dimensional cone structure formed by 3-interactions $\{((1,2),3),((2,1),3),((1,2),4),((2,1),4)\}$. $\beta_2^{L_1}=1$ corresponds to the connected component $\{1,2\}$. $\beta_2^{L_2}=1$ corresponds to the loop formed by $\{(1,2),(2,1)\}$. $\beta_3^{L_1}=1$ corresponds the connected component $\{(1,2),(2,1),3,4\}$. $\beta_3^{L_2}=1$ corresponds to the loop formed by 3-interactions $\{((1,2),3),((2,1),3),((1,2),4),((2,1),4)\}$.
	\begin{figure}[h]
		\centering
		\begin{minipage}{0.4\linewidth}
			\flushright
			\includegraphics[width=1\textwidth]{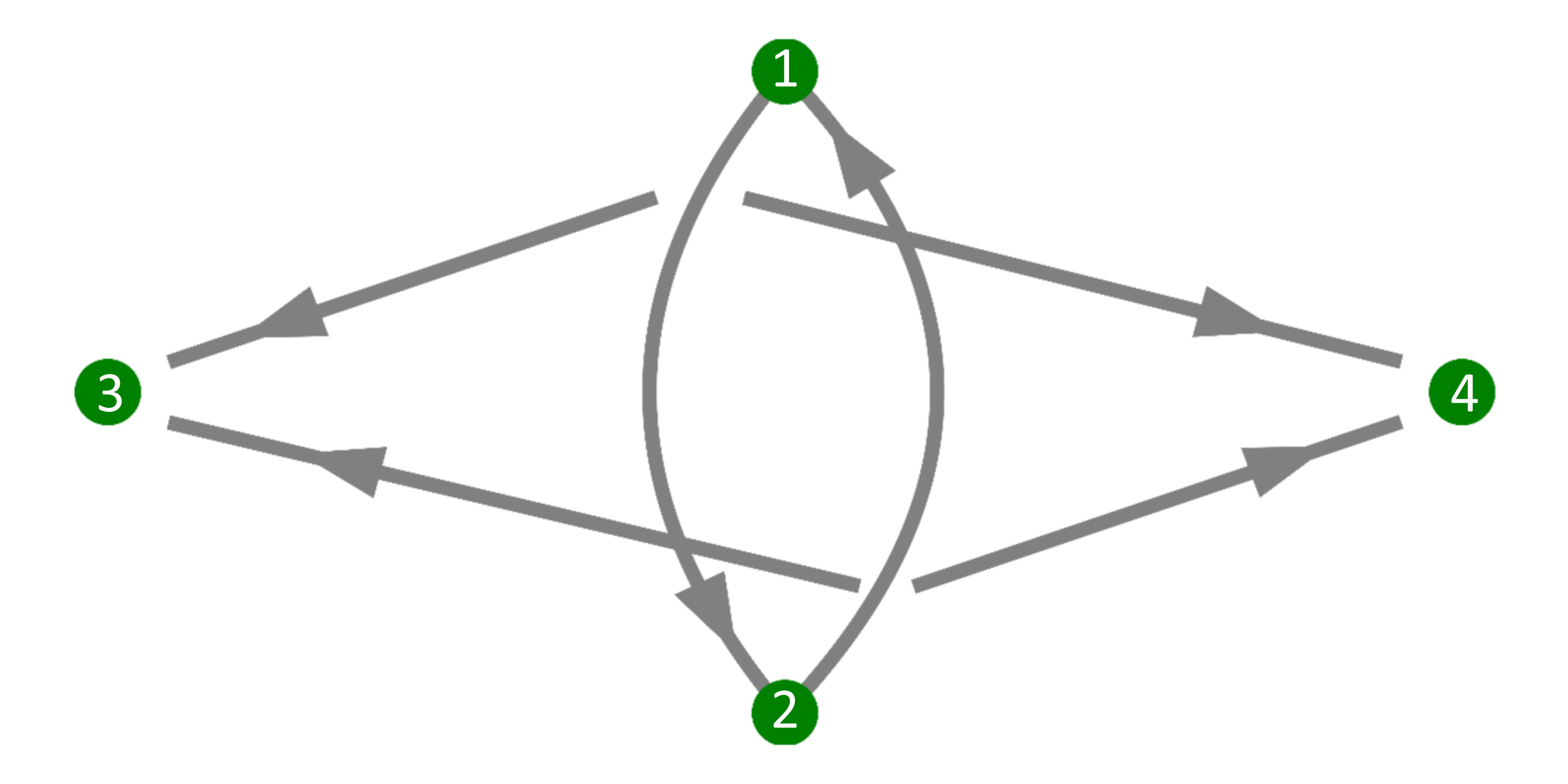}
			%\caption{An IntComplex.}
			%\label{fig:example2}
		\end{minipage}
		\hfill
		\begin{minipage}{0.5\linewidth}
			\flushleft
			\begin{tabular}{|c|c|c|c|c|}   
				\hline  &$\beta_1$ & $\beta_2$ & $\beta_3$ & $\beta_p(p>3)$ \\   
				\hline
				homology   &  3 & 0 & 1 & 0\\
				\hline
				Layer-Homology  &-& (1,1) & (1,1) & (0,0)\\
				\hline  	
			\end{tabular}   
			%\captionof{table}{Homology of the IntComplex in Figure \ref{fig:example2}. }
			%\label{table:example2}
		\end{minipage}
		\caption{An IntComplex and its homology.}
		\label{fig:4}
	\end{figure}
\end{example}

\begin{example}\label{exm:5}
	Consider the IntComplex in Figure \ref{fig:5}, $\mathcal{I}_1=V=\{1,2,3,4\}$, $\mathcal{I}_2=\{(1,2),(3,4)\}$, $\mathcal{I}_4=\{((1,2),(3,4)),((3,4),(1,2))\}$, $\mathcal{I}_p=\emptyset(p\ne 1,2,4)$. The homology are shown in the right table. $\beta_1=2$ corresponds to the connected components $\{1,2\}$ and $\{3,4\}$. $\beta_2^{L_1}=2$ corresponds to the connected components $\{1,2\}$ and $\{3,4\}$. $\beta_4^{L_1}=1$ corresponds to the connected component $\{(1,2),(3,4)\}$, $\beta_4^{L_2}=1$ corresponds to the loop formed by 4-interactions $\{((1,2),(3,4)),((3,4),(1,2))\}$.
	\begin{figure}[h]
		\centering
		\begin{minipage}{0.35\linewidth}
			\flushright
			\includegraphics[width=1\textwidth]{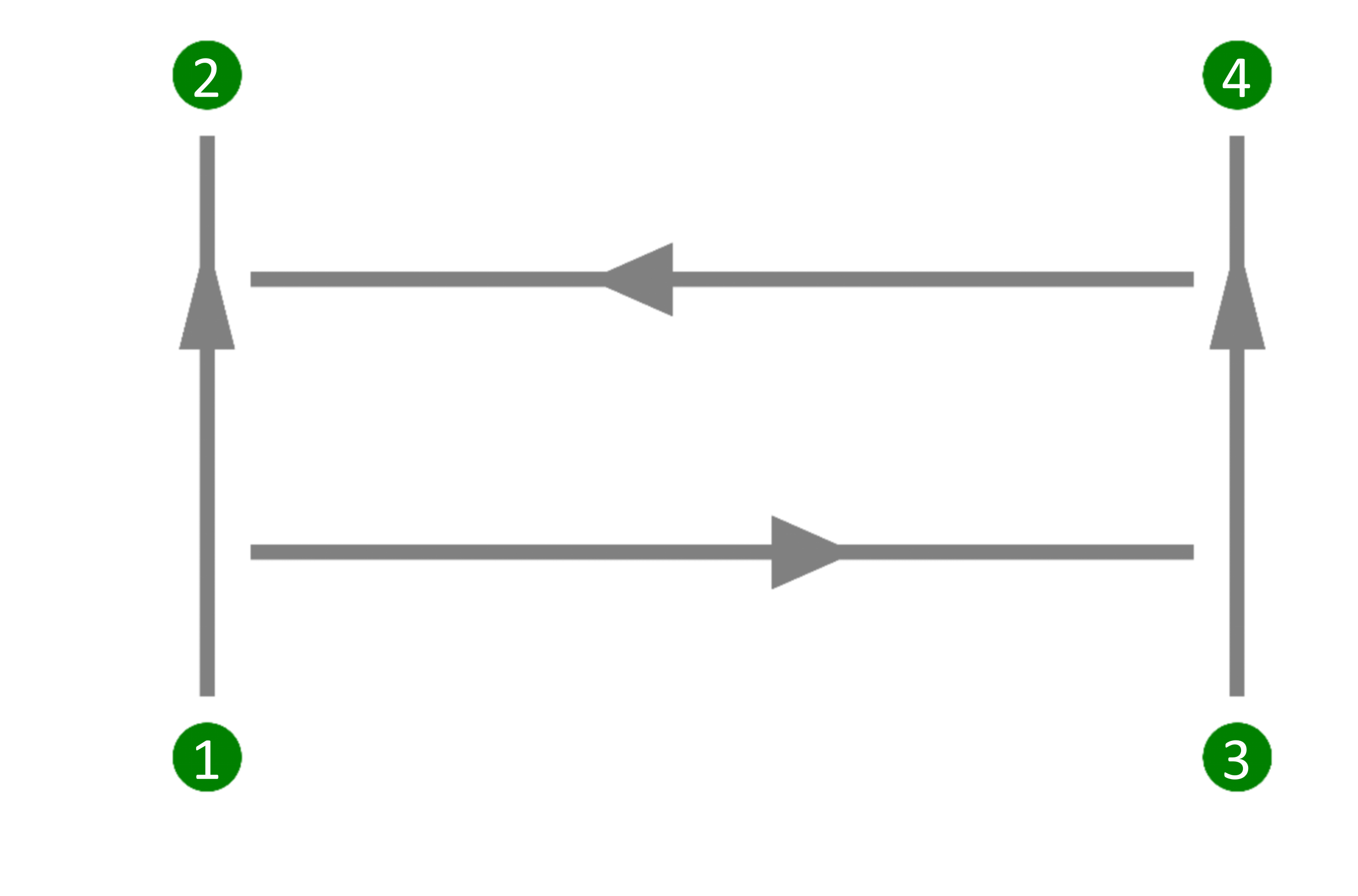}
			%\caption{An IntComplex.}
			%\label{fig:example2}
		\end{minipage}
		\hfill
		\begin{minipage}{0.58\linewidth}
			\flushleft
			\begin{tabular}{|c|c|c|c|c|c|}   
			\hline  &$\beta_1$ & $\beta_2$ & $\beta_3$ & $\beta_4$ & $\beta_p(p>4)$ \\   
			\hline
			homology   &  2 & 0 & 0 & 0 & 0\\
			\hline
			Layer-Homology  &-& (2,0) & (0,0) & (1,1) & (0,0)\\
			\hline  	
			\end{tabular}   
			%\captionof{table}{Homology of the IntComplex in Figure \ref{fig:example2}. }
			%\label{table:example2}
		\end{minipage}
		\caption{An IntComplex and its homology.}
		\label{fig:5}
	\end{figure}
\end{example}

\subsection{Multilayer-Homology}
The Multilayer-homology is for the structures composed of interactions in several layers, that is, a subset of the IntComplex that has interactions from several dimensions.
\begin{definition}[Multilayer-Homology]
	Given an IntComplex $\mathcal{I}$, consider its subset $S$ that contains interactions of several dimensions. Let $V_S^\mathcal{I}=\{\tau|\tau$ is the daughter of $\sigma,\sigma\in S\}$. $E_S^\mathcal{I}=\{\sigma=(\sigma_1,\sigma_2)|\sigma\in S\}$. Then $G^\mathcal{I}_S=(V^\mathcal{I}_S,E^\mathcal{I}_S)$ forms an IntComplex. We define the 2-th homology of $G^\mathcal{I}_S$ as the multilayer-homology of S.
	$$H(S)=H_2(G_S^\mathcal{I})$$
	Define the betti number of S as the dimension of the multilayer-homology, denoted by $\beta_S=dim(H_2(G_S^\mathcal{I}))$.
\end{definition}

The multilayer-homology of $S\subset \mathcal{I}$ reflects the loop structure formed by the interactions of $S$.  

\begin{example}\label{exm:6}
	Consider the IntComplex in Figure \ref{fig:6} {\bf A}, $\mathcal{I}_1=V=\{1,2,3,4,5,6,7,8,9\}$, $\mathcal{I}_2=\{(1,2),(4,5),(6,7),(8,9)\}$, $\mathcal{I}_3=\{((1,2),3),(3,(4,5)),(3,(6,7))\}$ $\mathcal{I}_4=\{((4,5),(8,9)),((6,7),(8,9))\}$, $\mathcal{I}_6=\{((1,2),((6,7),(8,9))),(((6,7),(8,9)),(4,5))\}$ $\mathcal{I}_p=\emptyset(p\ne 1,2,3,4,6)$. Assume $S=\mathcal{I}$, dimension of the multilayer-homology of $S$ is 2, which means the interactions of $S$ form two loops. One is formed by $((1,2),((6,7),(8,9)))$, $(((6,7),(8,9)),(4,5))$, $((1,2),3)$ and $(3,(4,5))$. The other is formed by $((6,7),(8,9))$, $((8,9),(4,5))$, $(3,(4,5))$ and $(3,(6,7))$, which are shown in the blue parts of Figure \ref{fig:6} {\bf B} and {\bf C}.
	\begin{figure}[h]
		\centering
		\includegraphics[width=1\textwidth]{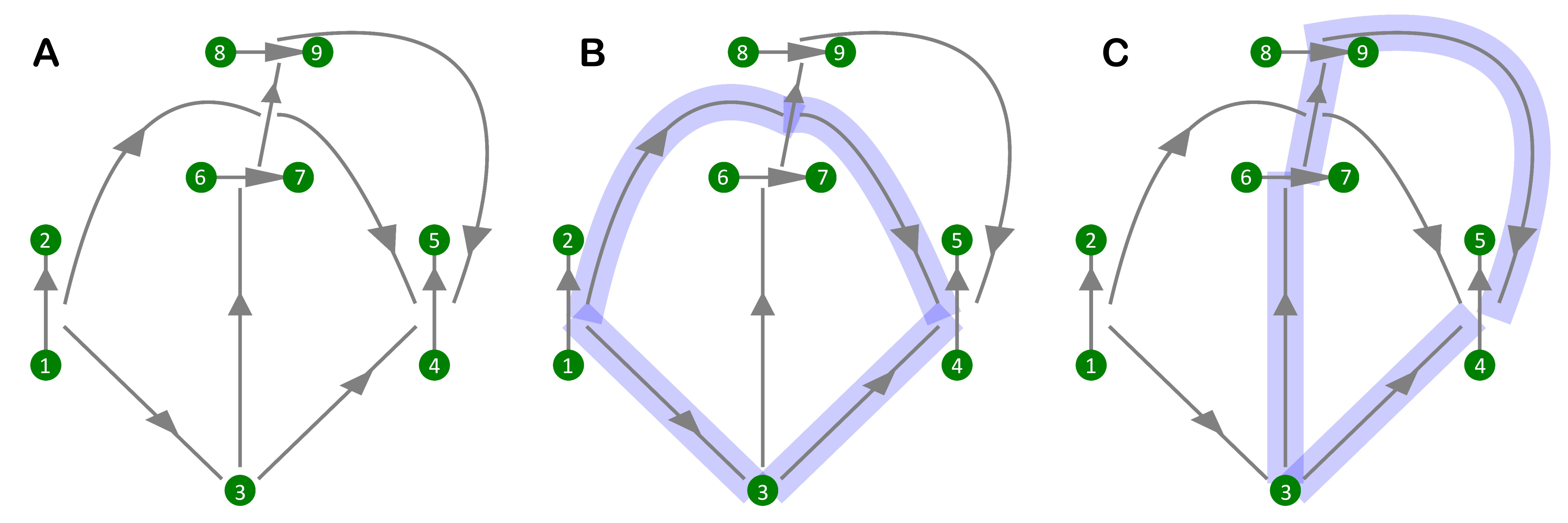}
		\caption{The subset $S$ of an IntComplex and its multilayer-homology. {\bf A}: An IntComplex over $V=\{1,2,3,4,5,6,7,8,9\}$. The blue loops in {\bf B} and {\bf C} corresponds to the two multilayer-homology classes of $S$ respectively.}
		\label{fig:6}
	\end{figure}
\end{example}

\subsection{Join}
\begin{definition}[Join of Interactions]
	Given two interactions $\sigma_p\in Int_p(V),\sigma_q\in Int_q(V)(1\leqslant p,q)$, define the join of $\sigma_p$ and $\sigma_q$ as the $(p+q)$-interaction $(\sigma_p,\sigma_q)\in Int_{p+q}(V).$
\end{definition}
Clearly, join of interactions is a bilinear operation that does not satisfy associative law. Assume $ \sigma_p=[v_{i_1}...v_{i_p},\{(l^a_s,r^a_s)\}_{s=1}^{p-1}]\in Int_p(V),\sigma_q=[v_{j_1}...v_{j_q},\{(l^b_t,r^b_t)\}_{t=1}^{q-1}]\in Int_q(V)$, then $(\sigma_p,\sigma_q)=[v_{i_1}...v_{i_p}v_{j_1}...v_{j_q},\{(l_s^a,r_s^a)\}_{s=1}^{p-1}\cup\{(l_t^b+p,r_t^b+p)\}_{t=1}^{q-1}\cup \{(1,p+q+1)\}]$.

\begin{prop}[Product rule]\label{prop:5}
	$\forall \sigma_p\in Int_p(V),\sigma_q\in Int_q(V)$, $\partial(\sigma_p,\sigma_q)=(\partial(\sigma_p),\sigma_q) + (-1)^p(\sigma_p,\partial(\sigma_q))$.
	\begin{proof}
		Assume $ \sigma_p=[v_{i_1}...v_{i_p},\{(l^a_s,r^a_s)\}_{s=1}^{p-1}],\sigma_q=[v_{j_1}...v_{j_q},\{(l^b_t,r^b_t)\}_{t=1}^{q-1}]$.
		\begin{equation}
			\begin{split}\nonumber
				\partial(\sigma_p)&=\partial([v_{i_1}...v_{i_p},\{(l_s^a,r_s^a)\}_{s=1}^{p-1}])\\ &=\sum_{k=1}^p(-1)^{k+1}[v_{i_1}...v_{i_{k-1}}v_{i_{k+1}}...v_{i_p},N_{\sigma_p}^k].\\
				(\partial(\sigma_p),\sigma_q)&=\sum_{k=1}^p(-1)^{k+1}[v_{i_1}...v_{i_{k-1}}v_{i_{k+1}}...v_{i_p}v_{j_1}...v_{j_q},N_{\sigma_p}^k\cup\{(l_t^b+p,r_t^b+p)\}_{t=1}^{q-1}\cup\{(1,p+q)\}].\\
				\partial(\sigma_q)&=\partial([v_{j_1}...v_{j_q},\{(l_t^b,r_t^b)\}_{t=1}^{q-1}])\\ &=\sum_{k=1}^q(-1)^{k+1}[v_{j_1}...v_{j_{k-1}}v_{j_{k+1}}...v_{j_q},N_{\sigma_q}^k]\\
				(-1)^p(\sigma_p,\partial(\sigma_q))&=\sum_{k=1}^q(-1)^{p+k}[v_{i_1}...v_{i_p}v_{j_1}...v_{j_{k_1}}v_{j_{k+1}}...v_{j_q},\{(l_s^a,r_s^a)\}_{s=1}^{p-1}\cup (N_{\sigma_q}^k+p)\cup\{(1,p+q)\}]
			\end{split}
		\end{equation}
		Here $N_{\sigma_q}^k+p=\{(a+p,b+p)|(a,b)\in N_{\sigma_q}^k\}$
		\begin{equation}\nonumber
			\begin{split}
				&\quad (\partial(\sigma_p),\sigma_q)+(-1)^p(\sigma_p,\partial(\sigma_q))\\
				&=\sum_{k=1}^p(-1)^{k+1}[v_{i_1}...v_{i_{k-1}}v_{i_{k+1}}...v_{i_p}v_{j_1}...v_{j_q},N_{\sigma_p}^k\cup\{(l_t^b+p,r_t^b+p)\}_{t=1}^{q-1}\cup\{(1,p+q)\}]\\
				&+\sum_{k=1}^q(-1)^{p+k}[v_{i_1}...v_{i_p}v_{j_1}...v_{j_{k_1}}v_{j_{k+1}}...v_{j_q},\{(l_s^a,r_s^a)\}_{s=1}^{p-1}\cup (N_{\sigma_q}^k+p)\cup\{(1,p+q)\}\cup\{(1,p+q)\}]\\
				&=\sum_{k=1}^p(-1)^{k+1}[v_{i_1}...v_{i_{k-1}}v_{i_{k+1}}...v_{i_p}v_{j_1}...v_{j_q},N_{\sigma_p}^k\cup\{(l_t^b+p,r_t^b+p)\}_{t=1}^{q-1}\cup\{(1,p+q)\}]\\
				&+\sum_{k=p+1}^{p+q}(-1)^{k}[v_{i_1}...v_{i_p}v_{j_1}...v_{j_{k_1}}v_{j_{k+1}}...v_{j_q},\{(l_s^a,r_s^a)\}_{s=1}^{p-1}\cup (N_{\sigma_q}^{k-p}+p)\cup\{(1,p+q)\}]\\
				&=\sum_{k=1}^{p}(-1)^{k+1}[v_{i_1}...v_{i_{k-1}}v_{i_{k+1}}...v_{i_p}v_{j_1}...v_{j_q},N_{(\sigma_p,\sigma_q)}^k]\\
				&+\sum_{k=p+1}^{p+q}(-1)^{k+1}[v_{i_1}...v_{i_p}v_{j_1}...v_{j_{k-1}}v_{j_{j+1}}...v_{j_q},N_{(\sigma_p,\sigma_q)}^k]\\
				&=\partial(\sigma_p,\sigma_q)
			\end{split}
		\end{equation}

	\end{proof}
\end{prop}

\subsection{Interaction map}
\begin{definition}[Interaction map]
Given two IntComplex $\mathcal{I},\mathcal{I}'$ over $V$ and $V'$ respectively, assume $f:V\rightarrow V'$ is a vertex map. For an interaction $\sigma=[v_1...v_k,\{(l_i,r_i)\}_{i=1}^{k-1}]\in \mathcal{I}$, let $f(\sigma)=[f(v_1)...f(v_k),\{(l_i,r_i)\}_{i=1}^{k-1}]$, if $f(\sigma)\in \mathcal{I}'$ for $\forall \sigma\in \mathcal{I}$, $f$ is called an interaction map from $\mathcal{I}$ to $\mathcal{I}'$.
\end{definition}
The interaction map $f:\mathcal{I}\rightarrow\mathcal{I}'$ sends an $n$-interaction of $\mathcal{I}$ to an $n$-interaction of $\mathcal{I}'$.  

\begin{prop}\label{prop:6}
	Given an interaction map $f:\mathcal{I}\rightarrow \mathcal{I}'$, $\forall \sigma=(\sigma_p,\sigma_q)\in \mathcal{I}$ $f(\sigma)=f(\sigma_p,\sigma_q)=(f(\sigma_p),f(\sigma_q))$
	\begin{proof}
		Assume $ \sigma_p=[v_{i_1}...v_{i_p},\{(l^a_s,r^a_s)\}_{s=1}^{p-1}],\sigma_q=[v_{j_1}...v_{j_q},\{(l^b_t,r^b_t)\}_{t=1}^{q-1}]$, we have 
		
		$(\sigma_p,\sigma_q)=[v_{i_1}...v_{i_p}v_{j_1}...v_{j_q},\{(l_s^a,r_s^a)\}_{s=1}^{p-1}\cup\{(l_t^b+p,r_t^b+p)\}_{t=1}^{q-1}\cup\{(1,p+q+1)\}]$.
		\begin{equation}\nonumber
			\begin{split}
				f(\sigma_p,\sigma_q)&=[f(v_{i_1})...f(v_{i_p})f(v_{j_1})...f(v_{j_q}),\{(l_s^a,r_s^a)\}_{s=1}^{p-1}\cup\{(l_t^b+p,r_t^b+p)\}_{t=1}^{q-1}\cup\{(1,p+q+1)\}]\\
				&=([f(v_{i_1})...f(v_{i_p}),\{(l_s^a,r_s^a)\}_{s=1}^{p-1}],[f(v_{j_1})...f(v_{j_q}),\{(l_t^b,r_t^b)\}_{t=1}^{q-1}])\\
				&=(f(\sigma_p),f(\sigma_q))
			\end{split}
		\end{equation}
	\end{proof}
\end{prop}
\begin{prop}\label{prop:7}
	Given an interaction map $f:\mathcal{I}\rightarrow \mathcal{I}'$, $f\partial=\partial f$.
	\begin{proof}
		$\forall \sigma_1\in \mathcal{I}_1,$ $f(\partial(\sigma_1))=0=\partial(f(\sigma_1))$.
		Assume $f(\partial(\sigma_p))=\partial(f(\sigma_p))$ for $\forall \sigma_p\in \mathcal{I}_p, \forall p<n$. Now we consider $n$-interactions, $\forall \sigma_n\in \mathcal{I}_n$, $\exists \sigma_p\in \mathcal{I}_p,\sigma_q\in\mathcal{I}_q,p,q<n$ such that $\sigma_n=(\sigma_p,\sigma_q)$.
		\begin{equation}
			\begin{split}
				\partial(f(\sigma_n))&=\partial(f(\sigma_p,\sigma_q))\\
				&=\partial(f(\sigma_p),f(\sigma_q))\\
				&=(\partial(f(\sigma_p),f(\sigma_q)+(-1)^p(f(\sigma_p),\partial(f(\sigma_q)))\\
				&=(f(\partial(\sigma_p)),f(\sigma_q)) + (-1)^p(f(\sigma_p),f(\partial(\sigma_q)))~~(by~ induction)\\
				&=f(\partial(\sigma_p),\sigma_q)+(-1)^pf(\sigma_p,\partial(\sigma_q))\\
				&=f(\partial(\sigma_p,\sigma_q))\\
				&=f(\partial(\sigma_n))
			\end{split}
		\end{equation}
	Hence the result follows.
	\end{proof}
\end{prop}
Consequently, $f:\mathcal{I}\rightarrow \mathcal{I}'$ induces a chain map $f_*:(\Lambda_*(\mathcal{I}),\partial_*)\rightarrow (\Lambda_*(\mathcal{I}'),\partial_*)$.
\begin{prop}\label{prop:8}
	Given an interaction map $f:\mathcal{I}\rightarrow \mathcal{I}'$, the map $f_*|_{\Omega_*(\mathcal{I})}$ provides a chain map 
	$$f_*|_{\Omega_*(\mathcal{I})}:\Omega_*(\mathcal{I})\rightarrow \Omega_*(\mathcal{I}')$$
	consequently, a homomorphism of homology groups
	$$f_\#:H_*(\mathcal{I})\rightarrow H_*(\mathcal{I}')$$
	\begin{proof}
		It suffices to prove that $f(\Omega_p(\mathcal{I}))\subset \Omega_p(\mathcal{I}')$. Firstly we have $f(\mathcal{A}_p(\mathcal{I}))\subset \mathcal{A}_p(\mathcal{I}')$. $\forall \sigma\in \Omega_p(\mathcal{I})$, we have $\sigma\in \mathcal{A}_p(\mathcal{I}),\partial(\sigma)\in \mathcal{A}_{p-1}(
		\mathcal{I})$, so $f(\sigma)\in \mathcal{A}_p(\mathcal{I}')$, $\partial(f(\sigma))=f(\partial(\sigma))\in \mathcal{A}_{p-1}(\mathcal{I}')$, which means $f(\sigma)\in \Omega_p(\mathcal{I}')$. Hence $f(\Omega_p(\mathcal{I}))\subset \Omega_p(\mathcal{I}')$. 
	\end{proof}
\end{prop}

%\begin{prop}
%	Any interaction map $f:\mathcal{I} \rightarrow \mathcal{I}'$ induces an interaction map 
%	$$f|_{G^\mathcal{I}_p}:G^\mathcal{I}_p\rightarrow G^{\mathcal{I}'}_p$$
%	and thus induces a homomorphism of Layer-Homology groups:
%	$$(f|_{G^\mathcal{I}_p})_\#:H^L_*(\mathcal{I})\rightarrow H^L_*(\mathcal{I}')$$
%	denoted by $f_\#^L$.
%	\begin{proof}
%		The interaction map sends an $n$-interaction to an $n$-interaction, so $f|_{G^\mathcal{I}_p}$ is an interaction map. The interaction map $f|_{G^\mathcal{I}_p}$ is the restriction of $f$ on $G_p^\mathcal{I}$, so it induces the homomorphism of Layer-Homology groups by proposition \ref{prop:8}. 
%	\end{proof}
%\end{prop}

\begin{prop}[Functoriality of Homology]\label{prop:14}
	Let $\mathcal{I},\mathcal{I}',\mathcal{I}''$ be three IntComplexes.
	\begin{enumerate}
		\item Let $id_\mathcal{I}:\mathcal{I}\rightarrow\mathcal{I}$ be the identity interaction map. Then $({id_\mathcal{I}}_\#)_p:H_p(\mathcal{I})\rightarrow H_p(\mathcal{I})$ is the identity linear map for each $p>0$.
		\item Let $f:\mathcal{I}\rightarrow \mathcal{I}'$, $g:\mathcal{I}'\rightarrow\mathcal{I}''$ be interaction maps. Then $((gf)_\#)_p=(g_\#)_p(f_\#)_p$ for each $p>0$.
	\end{enumerate}
	\begin{proof}
		For the first claim:
		$\forall \sigma_p=[v_1...v_p,\{(l_i,r_i)\}_{i=1}^{p-1}]\in \Omega_p(\mathcal{I})$,
		\begin{equation}\nonumber
			\begin{split}
				(id_{\mathcal{I}_*})_p(\sigma_p)&=[id_\mathcal{I}(v_1)...id_\mathcal{I}(v_p),\{(l_i,r_i)\}_{i=1}^{p-1}]\\
				&=[v_1...v_p,\{(l_i,r_i)\}_{i=1}^{p-1}]\\
				&=\sigma_p
			\end{split}
		\end{equation}  
	It follows that $(id_{\mathcal{I}_*})_p$ is the identity map on $\Omega_p(\mathcal{I})$, and thus $(id_{\mathcal{I}_\#})_p$ is the identity map on $H_p(\mathcal{I})$.
		
		For the second claim: $\forall \sigma_p=[v_1...v_p,\{(l_i,r_i)\}_{i=1}^{p-1}]\in \Omega_p(\mathcal{I})$,
		\begin{equation}\nonumber
			\begin{split}
				((gf)_*)_p(\sigma_p)&=[gf(v_1)...gf(v_p),\{(l_i,r_i)\}_{i=1}^{p-1}]\\
				&=[g(f(v_1))...g(f(v_p)),\{(l_i,r_i)\}_{i=1}^{p-1}]\\
				&=(g_*)_p([f(v_1)...f(v_p),\{(l_i,r_i)\}_{i=1}^{p-1}])\\
				&=(g_*)_p(f_*)_p([v_1...v_p,\{(l_i,r_i)\}_{i=1}^{p-1}])\\
				&=(g_*)_p(f_*)_p(\sigma_p)
			\end{split}
		\end{equation}
	The result follows. 
	\end{proof}
\end{prop}

%\begin{prop}[Functoriality of Layer-Homology]
%	Let $\mathcal{I},\mathcal{I}',\mathcal{I}''$ be three IntComplexes.
%	\begin{enumerate}
%		\item Let $id_\mathcal{I}:\mathcal{I}\rightarrow\mathcal{I}$ be the identity interaction map. Then $({id_\mathcal{I}}^L_\#)_p:H_p^L(\mathcal{I})\rightarrow H_p^L(\mathcal{I})$ is the identity linear map for each $p>0$.
%		\item Let $f:\mathcal{I}\rightarrow \mathcal{I}'$, $g:\mathcal{I}'\rightarrow\mathcal{I}''$ be interaction maps. Then $((gf)^L_\#)_p=(g^L_\#)_p(f^L_\#)_p$ for each $p>0$.
%	\end{enumerate}
%	\begin{proof}
%		The result follows from proposition \ref{prop:14}.
%	\end{proof}
%\end{prop}

\section{Persistent Homology of IntComplex}\label{section:persistent-homology}
\begin{figure}[h]
		\centering
		\includegraphics[width=1\textwidth]{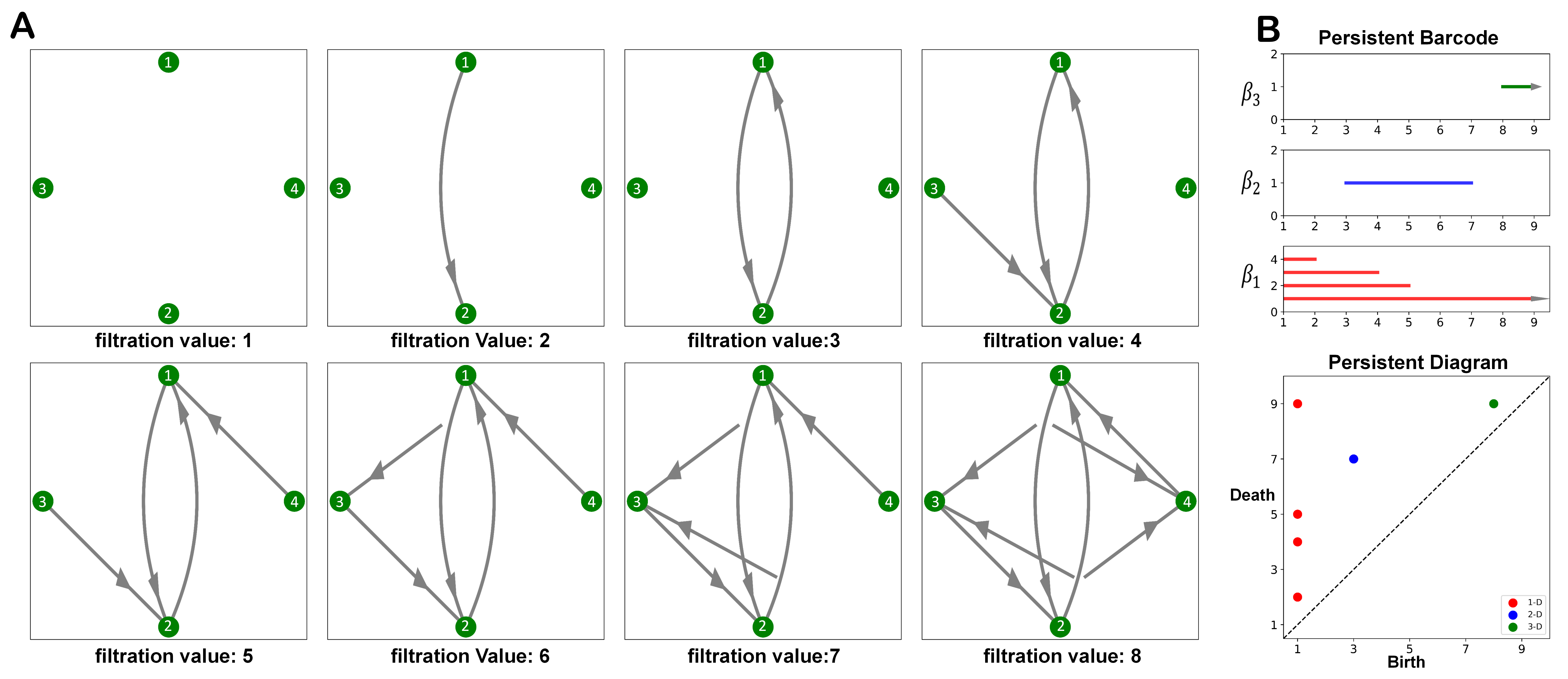}
		\caption{Illustration of the persistent homology of an IntComplex. {\bf A}: the filtration of an IntComplex. {\bf B}: the persistent barcode and peristent diagram representation of the persistent homology.}
		\label{fig:7}
\end{figure}

As a key theory of topological data analysis (TDA), persistent homology has been successfully applied to various fields. Here we give the persistent homology of IntComplex.
\subsection{Persistent homology}
\begin{definition}[IntComplex filtration]
	A filtration of the IntComplex $\mathcal{I}$ is a sequence of IntComplexes 
	$$\mathcal{I}(1)\subset\mathcal{I}(2)\subset...\subset\mathcal{I}(n)=\mathcal{I}$$
	where $\mathcal{I}(p)$ is a subcomplex of $\mathcal{I}(p+1)$, they are connected by inclusion map $(1\leqslant p<n).$ 
\end{definition}
The functoriality of homology enables us to obtain a persistent vector spaces from an IntComplex filtration, we give the definition as follows:
\begin{definition}[Persistent homology of IntComplex]
	Given an IntComplex filtration $\{\mathcal{I}(k)\xrightarrow{i_{k}}\mathcal{I}(k+1)\}_{0\leqslant k}$ of $\mathcal{I}$. Then for each $p>0$, the p-th persistent homology of the filtration is defined as the following persistent vector space:
	$$H_p^q(\mathcal{I})=Im\{H_p(\mathcal{I}(k))\xrightarrow{(i_{k,k+q})_\#}H_p(\mathcal{I}(k+q))\}_{q>0,k>0}$$
	Then define the p-th persistent diagram of the filtration $\mathcal{I}$ to be the persistent diagram of $H_p^q(\mathcal{I})$.
\end{definition}

%\begin{definition}[Persistent Layer-Homology of IntComplex]
%	Given an IntComplex filtration $\{\mathcal{I}(k)\xrightarrow{i_{k}}\mathcal{I}(k+1)\}_{0\leqslant k}$ of $\mathcal{I}$. Then for each $p>0$, the p-th persistent Layer-Homology of the filtration is defined as the following two persistent vector space:
%	$$H_p^{L_1,q}(\mathcal{I})=Im\{H_p^{L_1}(\mathcal{I}(k))\xrightarrow{(i_{k,k+q})_\#}H_p^{L_1}(\mathcal{I}(k+q))\}_{q>0,k>0}$$
	
%	$$H_p^{L_2,q}(\mathcal{I})=Im\{H_p^{L_2}(\mathcal{I}(k))\xrightarrow{(i_{k,k+q})_\#}H_p^{L_2}(\mathcal{I}(k+q))\}_{q>0,k>0}$$
%	Then define the p-th persistent layer-diagram of the filtration $\mathcal{I}$ to be the persistent diagram of $H_p^{L_1,q}(\mathcal{I})$ and $H_p^{L_2,q}(\mathcal{I})$.
%\end{definition}
\begin{definition}[Weighted InComplex]\label{weighted-complex}
	A weighted IntComplex is a pair $(\mathcal{I},f)$ where $\mathcal{I}$ is an IntComplex $\mathcal{I}$ and $f$ is a function $f:\mathcal{I}\rightarrow R$, denoted by $\mathcal{I}^f$. let $\mathcal{I}(r)=\{\sigma|f(\sigma)\leqslant r\}$, then we get an IntComplex filtration $\{\mathcal{I}(r)\xrightarrow{i}\mathcal{I}(s)\}_{r<s}$. We denote the persistent homology and persistent diagram of this filtration by $H(\mathcal{I},f)$ and $\mathcal{D}(\mathcal{I}^f)$ respectively.
\end{definition}
%\begin{prop}\label{prop:15}
	%Given a weighted IntComplex $(\mathcal{I},f)$, let $\mathcal{I}(r)=\{\sigma|f(\sigma)\leqslant r\}$. Then we get an IntComplex filtration $\{\mathcal{I}(r)\xrightarrow{i}\mathcal{I}(s)\}_{r<s}$. We denote the persistent diagram of this filtration by $\mathcal{D}(\mathcal{I},f)$
	%\begin{proof}
	%	This is an obvious fact.
%	\end{proof}
%\end{prop}

\begin{example}\label{exp:7}
	Consider a weighted IntComplex $\mathcal{I}$ over $V=\{1,2,3,4\}$, $\mathcal{I}=\{[1,1]$, $[2,1]$, $[3,1]$, $[4,1]$, $[(1,2),2]$, $[(2,3),3]$, $[(3,2),4]$, $[(4,1),5]$, $[((1,2),3),6]$, $[((2,1),3),7]$, $[((1,2),4),8]$, $[((2,1),4),8]\}$ where $[\sigma,v]$ represents an interaction $\sigma$ with weight $v$, denoted by $w(\sigma)=v$. Let $\mathcal{I}(k)=\{\sigma|w(\sigma)\leqslant k\}$. Then we get a filtration of $\mathcal{I}$,
	$$\mathcal{I}(1)\subset\mathcal{I}(2)\subset\mathcal{I}(3)\subset\mathcal{I}(4)\subset\mathcal{I}(5)\subset\mathcal{I}(6)\subset\mathcal{I}(7)\subset\mathcal{I}(8)$$
	which is shown in Figure \ref{fig:7} {\bf A}. The persistent homology is computed and represented by persistent barcode and persistent diagram as in Figure \ref{fig:7} {\bf B}.
\end{example}

\begin{figure}[ht]
		\centering
		\includegraphics[width=0.95\textwidth]{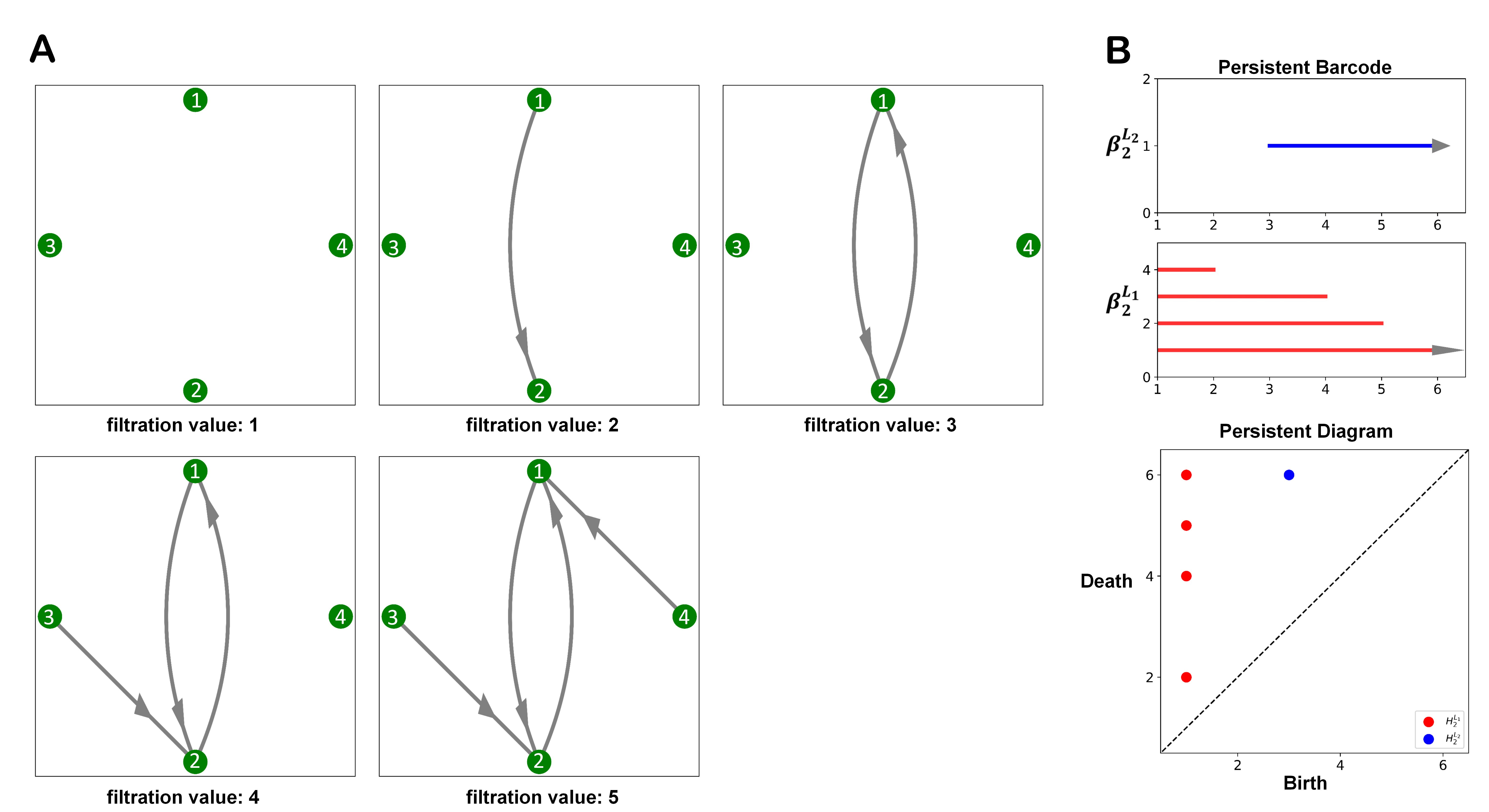}
		\caption{Illustration of the persistent $p$-layer-homology of an weighted IntComplex. {\bf A}: the induced filtration for $p=2$. {\bf B}: the persistent barcode and peristent diagram representation of the persistent 2-layer-homology.}
		\label{fig:8}
\end{figure}

\begin{figure}[ht]
		\centering
		\includegraphics[width=0.95\textwidth]{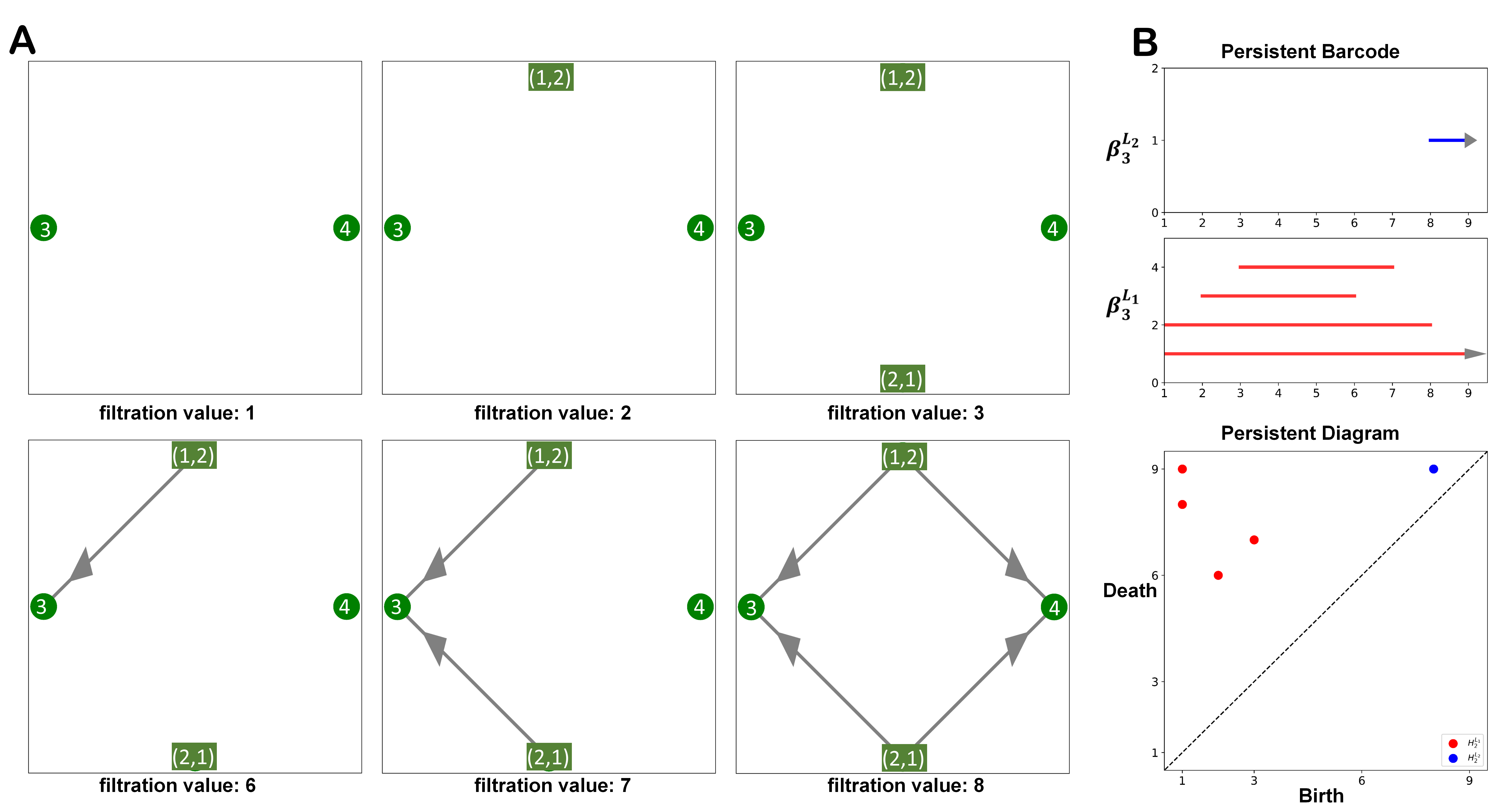}
		\caption{Illustration of the persistent $p$-layer-homology of an weighted IntComplex. {\bf A}: the induced filtration for $p=3$. {\bf B}: the persistent barcode and peristent diagram representation of the persistent 3-layer-homology.}
		\label{fig:9}
	\end{figure}

\begin{figure}[ht]
		\centering
		\includegraphics[width=0.95\textwidth]{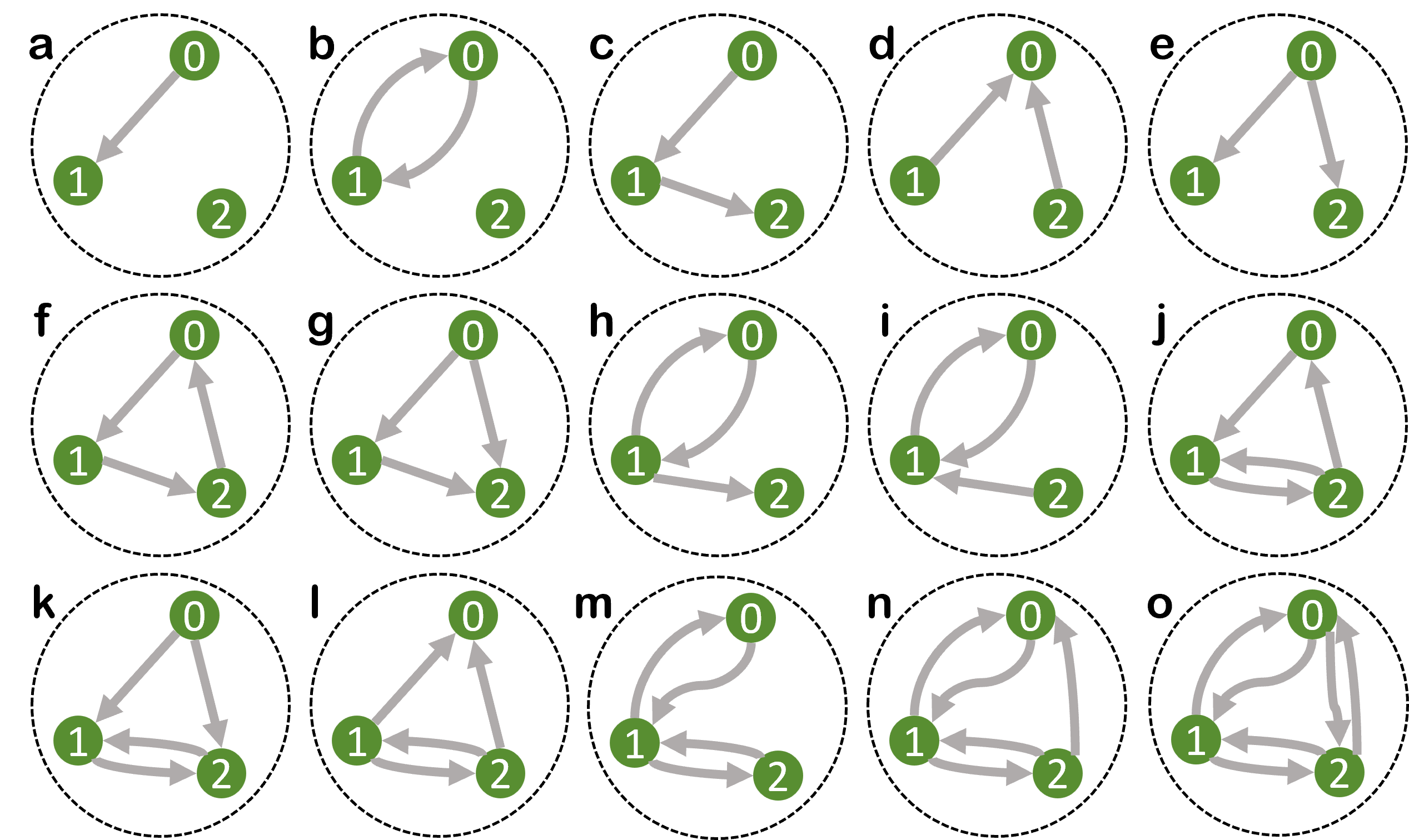}
		\caption{Fifteen types of digraphs with at least one arrow over three vertices.}
		\label{fig:12}
\end{figure}

\subsection{Stability of persistent homology}
The persistent diagram of an IntComplex filtration is a multi-set of pairs in which each pair $(b_i,d_i)$ represents a homology class that appears at filtration $b_i$ and disappears at filtration $d_i$. Note that $d_i$ can be $\infty$ and we define $|(+\infty)-(+\infty)|=|(-\infty)-(-\infty)|=0$. For any two pairs $p_i=(b_i,d_i),p_j=(b_j,d_j)$, define $d^\infty(p_i,p_j)=\max\{|b_i-b_j|,|d_i-d_j|\}$ and $d^\infty(p_i)=(d_i-b_i)/2$. 
\begin{definition}[Partial matching]
    A partial matching between multi-sets $A$ and $B$ is a collection of pairs $\psi=\{(\alpha,\beta)\in A\times B\}$ where $\alpha$ and $\beta$ can occur in at most one pair. If $(\alpha,\beta)\in \psi$, we say that $\alpha$ is matched with $\beta$, otherwise if $\alpha$ does not belong to any pair in $\psi$, we say $\alpha$ is unmatched.
\end{definition}
\begin{definition}[$\delta$-matching]
    A partial matching $\psi$ between A and B is called a $\delta$-matching if
    \begin{itemize}
        \item $\forall (\alpha,\beta)\in \psi$, $d^\infty(\alpha,\beta)\leqslant \delta$
        \item if $\alpha\in A\cup B$ is unmatched, $d^\infty(\alpha)\leqslant\delta$
    \end{itemize}
\end{definition}
\begin{definition}[Bottleneck distance]
    The bottleneck distance between two multi-sets $A$ and $B$ is defined as 
    $$d_B(A,B)=\inf\{\delta|there~exists~a~\delta\text{-}matching~between~A~and~B\}$$
\end{definition}
As with any category, a homomorphism of degree $\epsilon$ between two persistence modules $\mathbb{U}=\{U_s\xrightarrow{u_{s}^t}{U_t}\}$ and $\mathbb{V}=\{V_s\xrightarrow{v_{s}^t}V_t\}$ is a collection $\Phi$ of module morphisms 
$\phi_t:U_t\rightarrow V_{t+\epsilon}$ such that $\phi_tu_{s}^t=v_{s+\epsilon}^{t+\epsilon}\phi_s$ for all $s\leqslant t$. Let ${\rm Hom}^\epsilon(\mathbb{U},\mathbb{V})$ be the collection of all homomorphism of degree $\epsilon$ from $\mathbb{U}$ to $\mathbb{V}$.
\begin{definition}[$\epsilon$-interleaved]
    Two persistence modules $\mathbb{U}$ and $\mathbb{V}$ are called $\epsilon$-interleaved if there are maps
    $$\Phi\in{\rm Hom}^\epsilon(\mathbb{U},\mathbb{V}),\Psi\in {\rm Hom}^\epsilon(\mathbb{V},\mathbb{U})$$
    such that
    $$\Phi\Psi=1_\mathbb{V}^{2\epsilon},\Psi\Phi=1_\mathbb{U}^{2\epsilon}$$
\end{definition}

\begin{definition}[Interleaving distance]
    The interleaving distance between two persistence modules $\mathbb{U}$ and $\mathbb{V}$ is defined as
    $$d_I(\mathbb{U},\mathbb{V})=\inf\{\epsilon|\mathbb{U}~and~\mathbb{V}~are~ \epsilon\text{-}interleaved\}$$
\end{definition}

Let $f,g$ be two real-valued functions defined on an finite IntComplex $\mathcal{I}$, $||f-g||_\infty=\sup\limits_{\sigma\in\mathcal{I}}{|f(\sigma)-g(\sigma)|}$, $\mathcal{D}(\mathcal{I}^f)$ and $\mathcal{D}(\mathcal{I}^g)$ be the persistence diagrams of $\mathcal{I}$ induced by $f$ and $g$ respectively. We have the following result:
\begin{thm}[Stability of persistent homology]
	$$d_B(\mathcal{D}(\mathcal{I}^f),\mathcal{D}(\mathcal{I}^g))\leqslant||f-g||_\infty$$.
 \begin{proof}
     From the algebraic stability, we have
     $$d_I(H(\mathcal{I},f),H(\mathcal{I},g))=d_I(\mathcal{D}(\mathcal{I}^f),\mathcal{D}(\mathcal{I}^g))$$
     Considering the category $\mathcal{C}$ of filtered IntComplex in which objects are filtered IntComplexs and morphisms are the interaction maps between them. $\mathcal{I}^f,\mathcal{I}^g\in\mathcal{C}$. Let $\epsilon=||f-g||_\infty$, we prove that $\mathcal{I}^f$ and $\mathcal{I}^g$ are $\epsilon$-interleaved in $\mathcal{C}$.

     We have the inclusions of IntComplexs $\mathcal{I}^f_a\hookrightarrow\mathcal{I}^f_{a+\epsilon}$, $\mathcal{I}^g_a\hookrightarrow\mathcal{I}^g_{a+\epsilon}$. Note that $\mathcal{I}^f_a$ and $\mathcal{I}^g_a$ are subcomplexes of $\mathcal{I}^g_{a+\epsilon}$ and $\mathcal{I}^f_{a+\epsilon}$ respectively. So we have the homomorphisms 
     $$\Phi:\mathcal{I}^f_\star\hookrightarrow \mathcal{I}^g_{\star+\epsilon},\Psi:\mathcal{I}^g_\star\hookrightarrow\mathcal{I}^f_{\star+\epsilon}$$
     All the morphisms are inclusions, so we have 
     $$\Psi\Phi=1_{\mathcal{I}^f_\star}^{2\epsilon},\Phi\Psi=1_{\mathcal{I}^g_\star}^{2\epsilon}$$
     Consequently, $\mathcal{I}^f$ and $\mathcal{I}^g$ are $\epsilon$-interleaved, which means $$d_I(\mathcal{I}^f,\mathcal{I}^g)\leqslant\epsilon$$
     by the functoriality of IntComplex homology, we have
     $$d_I(H(\mathcal{I}^f),H(\mathcal{I}^g))\leqslant d_I(\mathcal{I}^f,\mathcal{I}^g)\leqslant \epsilon$$
     So
     $$d_I(\mathcal{D}(\mathcal{I}^f)\mathcal{D}(\mathcal{I}^g))\leqslant\epsilon$$
     from the result in, we have
     $$d_B(\mathcal{D}(\mathcal{I}^f),\mathcal{D}(\mathcal{I}^g)\leqslant||f-g||_\infty$$
 \end{proof}
\end{thm}

\begin{definition}[Persistent $p$-layer-homology of weighted IntComplex]
	Given a weighted IntComplex $\hat{\mathcal{I}}=(\mathcal{I},f)$, consider $G_p^\mathcal{I}$, we get a weighted IntComplex $\mathcal{I}'=(G_p^\mathcal{I},f|_{G_p^\mathcal{I}})$. Then we can get a IntComplex filtration from $\mathcal{I}'$ by definition \ref{weighted-complex}, the persistent homology of this filtration is defined as the persistent $p$-layer homology of $\hat{\mathcal{I}}$. 
\end{definition}

\begin{example}
	For the weighted IntComplex $\mathcal{I}$ in Example \ref{exp:7}, we consider the persistent $p$-layer homology, the induced filtration for $p=2$ is shown in  Figure \ref{fig:8} {\bf A}. And the associated persistent $2$-layer homology is illustrated by persistent barcode and persistent diagram in Figure \ref{fig:8} {\bf B}.
 
	The induced filtration for $p=3$ is shown in Figure \ref{fig:9} {\bf A}. And the associated persistent $3$-layer homology is illustrated by persistent barcode and persistent diagram in Figure \ref{fig:9} {\bf B}.

	Note that for $p=3$, the 1-interaction set is $\{3,4,(1,2),(2,1)\}$, 2-interaction set is all the 3-interactions of $\mathcal{I}$. 
\end{example}

\section{Experimental analysis}\label{section:experiment}
In this section, we give experimental test to show the importance of high-order interactions for network analysis.

Considering the digraphs without self-loops over three vertices. There are 15 types of such digraphs with at least one arrow up to graph isomorphism, which is shown in Figure \ref{fig:12}.

\begin{figure}[h]
		\centering
		\includegraphics[width=0.95\textwidth]{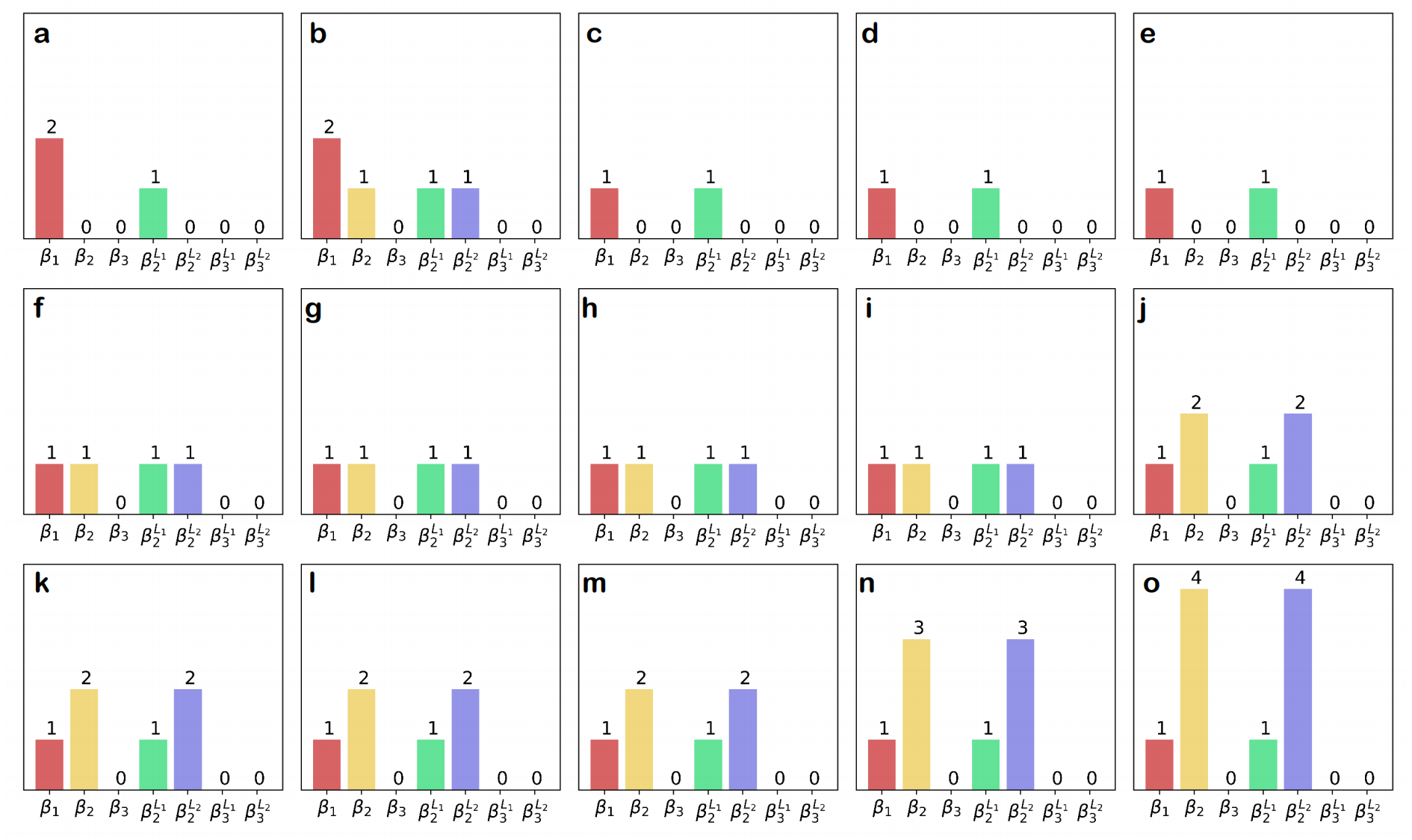}
		\caption{Homology, 2-layer homology and 3-layer homology of the 15 digraphs in Figure \ref{fig:12}. It can be seen that some non-isomorphic digraphs share the same homology information.}
		\label{fig:13}
\end{figure}

\begin{figure}[ht]
		\centering
		\includegraphics[width=0.95\textwidth]{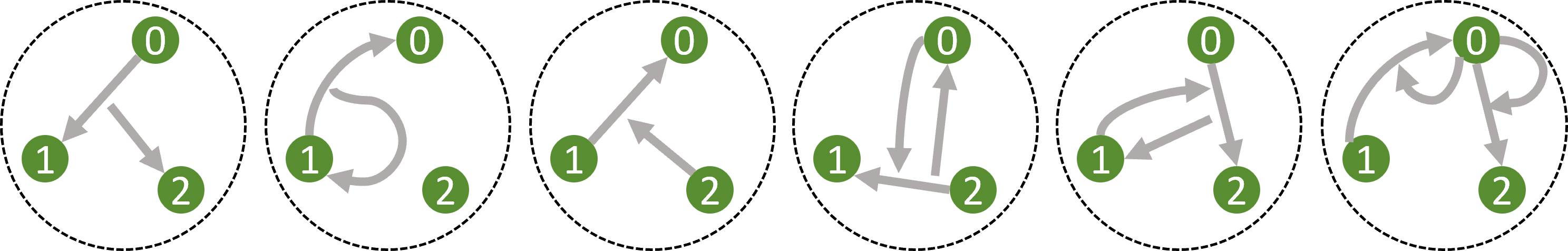}
		\caption{The 3-interactions that are added to the 15 digraphs in Figure \ref{fig:13}.}
		\label{fig:15}
\end{figure}
\begin{figure}[ht]
		\centering
		\includegraphics[width=0.95\textwidth]{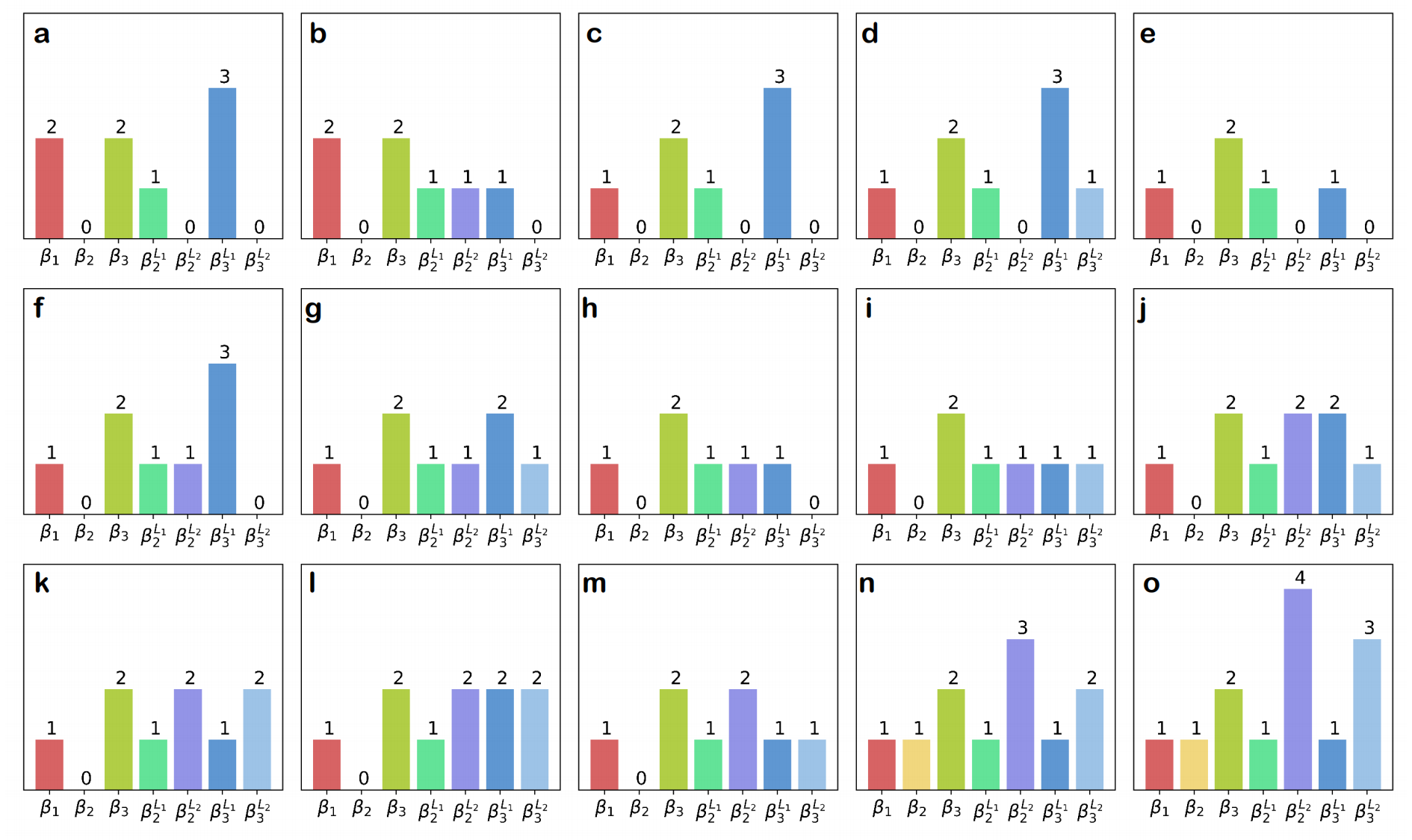}
		\caption{Homology, 2-layer homology and 3-layer homology of the IntComplexes derived from the 15 digraphs by adding the 3-interactions in Figure \ref{fig:15}. It can be seen that the homology information can clearly differentiate all 15 IntComplexes.}
		\label{fig:14}
\end{figure}
Among these digraphs, there is one digraph with one arrow ({\bf a}), four digraphs with two arrows ({\bf b,c,d,e}), four digraphs with three arrows ({\bf f,g,h,i}), four digraphs with four arrows ({\bf j,k,l,m}), one digraph with five arrows ({\bf n}) and one digraph with six arrows ({\bf o}).

Each digraph $G=(V,E)$ can be seen as an IntComplex $\mathcal{I}$ where $\mathcal{I}_0=V, \mathcal{I}_1=E$ and $\mathcal{I}_p=\emptyset(p>1)$. We compute the homology, 2-layer homology and 3-layer homology of these 15 digraphs. The results are shown in Figure \ref{fig:13}. It can be seen that the four digraphs with three arrows ({\bf f,g,h,i}) have the same homology information, the four digraphs with four arrows ({\bf j,k,l,m}) have the same homology information, and three digraphs (\textbf{c,d,e}) with two arrows have the same homology information, which means the homology information cannot differentiate them. 

As IntComplexes, these digraphs only have 1-interaction of independent identities and 2-interaction of pair-wise interactions. We add the following 3-interactions into all 15 IntComplexes to get 15 new IntComplexes with high-order interactions: ((0,1),2), ((1,0),1), (2,(1,0)), (0,(2,1)), ((2,1),0), (1,(0,2)), ((0,2),1), (0,(1,0)), (0,(0,2)). Figure \ref{fig:15} shows all these 3-interactions.

Then, we compute the homology, 2-layer homology and 3-layer homology of the new IntComplexes, the results are shown in Figure \ref{fig:14}.

It can be seen that there are no two IntComplexes have the same homology information for all 15 IntComplexes. In a summary, before adding the 3-interactions, the fifteen digraphs cannot be differentiated by the homology information, while after adding some 3-interactions, all 15 digraphs can be clearly differentiated by the homology information.

%\begin{table}[ht]
%	\caption{Homology and layer homology of 15 digraphs over three vertices.}
%	\label{table:jarvis-result}
%	\centering
%	\begin{tabular}{l|ccccccccccccccc}
%		\hline
%		  & a  & b & c & d & e & f & g & h & i & j & k& l & m & n & o \\
%		\hline
%		$\beta_1$  & 2 &2&1&1&1&1&1&1&1&1&1&1&1&1&1 \\
%		\hline
%		$\beta_2$ &  0 &1&0&0&0&1&1&1&1&2&2&2&2&3&4  \\
%		\hline	
 %       $\beta_3$ &  0 &0&0&0&0&0&0&0&0&0&0&0&0&0&0  \\
  %      \hline
%        $\beta_2^{L_1}$ & 1&1&1&1&1&1&1&1&1&1&1&1&1&1&1 \\
%        \hline
%        $\beta_2^{L_2}$ & 0&1&0&0&0&1&1&1&1&2&2&2&2&3&4\\
%       \hline
%        $\beta_3^{L_1}$ & 0&0&0&0&0&0&0&0&0&0&0&0&0&0&0\\
%        \hline
%        $\beta_3^{L_2}$ & 0&0&0&0&0&0&0&0&0&0&0&0&0&0&0\\
%        \hline
%	\end{tabular}
%\end{table}

\section{Conclusion}\label{section:conclusion}
In the work, we present IntComplex as a novel model for high-order interactions that existing graph, simplicial complex, and hypergraph models cannot represent. We introduce the homology theory, including the standard homology for adjacent order interactions, layer-homology for a specific order interactions and multi-layer homology for interactions across several orders, to give quantitative characterization of the IntComplex, enabling a detailed dissection of the topological structure of high-order interactions. Further, we introduce the persistent homology by considering the filtration process and give the stability result to ensure its robustness. IntComplex provides a foundational framework for analyzing the topological properties of high-order interactions, offering significant potential for advancing complex network analysis.

%\bibliography{refs}
%\bibliographystyle{unsrt}

\end{document}